\documentclass[leqno,12pt]{article}

\usepackage{amssymb,amscd}

\usepackage{amsthm}

\theoremstyle{plain} 

\newtheorem{theorem}{\sc Theorem}[section]

\newtheorem{lemma}[theorem]{\sc Lemma}

\newtheorem{corollary}[theorem]{\sc Corollary}

\newtheorem{proposition}[theorem]{\sc Proposition}

\theoremstyle{definition}

\newtheorem{remark}[theorem]{\sc Remark}

\newtheorem{example}[theorem]{\sc Example}

\newtheorem{remarks}[theorem]{\sc Remarks}

\newcommand\bF{{\mathbb F}}
\newcommand\bG{{\mathbb G}}

\newcommand\bP{{\mathbb P}}

\newcommand\bZ{{\mathbb Z}}

\newcommand\cA{{\mathcal A}}

\newcommand\cE{{\mathcal E}}

\newcommand\cG{{\mathcal G}}
\newcommand\cH{{\mathcal H}}

\newcommand\cL{{\mathcal L}}
\newcommand\cM{{\mathcal M}}

\newcommand\cO{{\mathcal O}}
\newcommand\cP{{\mathcal P}}

\newcommand\cX{{\mathcal X}}

\newcommand\fpgl{\mathfrak{pgl}}

\newcommand\tg{\tilde{g}}

\newcommand\tx{\tilde{x}}
\newcommand\ty{\tilde{y}}

\newcommand\tH{\widetilde{H}}
\newcommand\tHs{\widetilde{H_s}}
\newcommand\tHso{\widetilde{H_s^0}}
\newcommand\tHu{\widetilde{H_u}}

\newcommand\trho{\tilde{\rho}}

\newcommand\wX{\widehat{X}}

\newcommand\ad{{\rm ad}}

\newcommand\gp{{\rm gp}}

\newcommand\id{{\rm id}}

\newcommand\op{{\rm op}}

\newcommand\red{{\rm red}}

\newcommand\Aut{{\rm Aut}}
\newcommand\Br{{\rm Br}}
\newcommand\Div{{\rm Div}}
\newcommand\End{{\rm End}}

\newcommand\GL{{\rm GL}}
\newcommand\GO{{\rm GO}}
\newcommand\Gr{{\rm Gr}}
\newcommand\Hilb{{\rm Hilb}}
\newcommand\Hom{{\rm Hom}}
\newcommand\Int{{\rm Int}}

\newcommand\Lie{{\rm Lie}}
\newcommand\M{{\rm M}}
\newcommand\NS{{\rm NS}}
\renewcommand\O{{\rm O}}
\renewcommand\P{{\rm P}}
\newcommand\PGL{{\rm PGL}}
\newcommand\Pic{{\rm Pic}}
\newcommand\PO{{\rm PO}}
\newcommand\PSO{{\rm PSO}}
\newcommand\PSp{{\rm PSp}}
\newcommand\Proj{{\rm Proj}}
\newcommand\SL{{\rm SL}}
\newcommand\Sp{{\rm Sp}}
\newcommand\SO{{\rm SO}}



\title{Homogeneous projective bundles over abelian varieties}

\author{Michel Brion}

\date{}

\begin{document}

\maketitle

\footnote{%2010 MSC numbers
}

\begin{abstract}
We consider projective bundles (or Brauer-Severi varieties) 
over an abelian variety which are homogeneous, i.e., 
invariant under translation. We describe the structure 
of these bundles in terms of projective representations 
of commutative group schemes; the irreducible bundles 
correspond to Heisenberg groups and their standard 
representations. Our results extend those of Mukai on 
semi-homogeneous vector bundles, and yield a geometric 
view of the Brauer group of abelian varieties.
\end{abstract}

\section{Introduction}
\label{sec:introduction}

The main objects of this article are the projective bundles
(or Brauer-Severi varieties) over an abelian variety $X$ which
are homogeneous, i.e., isomorphic to their pull-backs under 
all translations. Among these bundles, the projectivizations 
of vector bundles are well understood by work of Mukai 
(see \cite{Muk78}).
Indeed, the vector bundles with homogeneous projectivization 
are exactly the semi-homogeneous vector bundles of [loc.~cit.]. 
Those that are simple (i.e., their global endomorphisms are just 
scalars) admit several remarkable characterizations; for example, 
they are all obtained as direct images of line bundles under 
isogenies. Moreover, every indecomposable semi-homogeneous vector 
bundle is the tensor product of a unipotent bundle and of a simple 
semi-homogeneous bundle.

\medskip

In this article, we obtain somewhat similar statements for the 
structure of homogeneous projective bundles. We build on 
the results of our earlier paper \cite{Br12} about homogeneous 
principal bundles under an arbitrary algebraic group; here we 
consider of course the projective linear group $\PGL_n$.
In loose terms, the approach of [loc.~cit.] reduces the
classification of homogeneous bundles to that of
commutative subgroup schemes of $\PGL_n$. The latter, 
carried out in Section \ref{sec:str}, is based on the classical 
construction of Heisenberg groups and their irreducible 
representations.

\medskip

In Section \ref{sec:irr}, we introduce a notion of irreducibility 
for homogeneous projective bundles, which is equivalent to the 
group scheme of bundle automorphisms being finite. 
(The projectivization of a semi-homogeneous vector bundle 
$E$ is irreducible if and only if $E$ is simple). 
We characterize those projective bundles that are 
homogeneous and irreducible, by the vanishing of all the 
cohomology groups of their adjoint vector bundle 
(Proposition \ref{prop:homirr}). Also, we show that
the homogeneous irreducible bundles are classified by 
the pairs $(H,e)$, where $H$ is a finite subgroup of the 
dual abelian variety, and $e: H \times H \to \bG_m$ a 
non-degenerate alternating bilinear pairing (Proposition 
\ref{prop:irraut}). Finally, we obtain a characterization
of those homogeneous projective bundles that are 
projectivizations of vector bundles, first in the irreducible
case (Proposition \ref{prop:semi}; it states in loose terms
that the pairing $e$ originates from a line bundle on $X$) 
and then in the general case (Theorem \ref{thm:obs}).

\medskip

The irreducible homogeneous projective bundles over an elliptic
curve are exactly the projectivizations of indecomposable
vector bundles with coprime rank and degree, as follows 
from classical work of Atiyah (see \cite{At57}). But any abelian 
variety $X$ of dimension at least $2$ admits many homogeneous 
projective bundles that are not projectivizations of vector 
bundles. In fact, any class in the Brauer group $\Br(X)$ is
represented by a homogeneous bundle (as shown by Elencwajg 
and Narasimhan in the setting of complex tori, see 
\cite[Thm.~1]{EN83}). Also, our approach yields a geometric 
view of a description of $\Br(X)$  due to Berkovich 
(see \cite{Be72}); this is developed in Remark \ref{rem:ber}. 

\medskip
 
Spaces of algebraically equivalent effective divisors on an 
arbitrary projective variety afford geometric examples of projective 
bundles. These spaces are investigated in Section \ref{sec:exa}
for abelian varieties and curves of genus $g \geq 2$; they turn out 
to be homogeneous in the former case, but not in the latter. 

\medskip

In the final Section \ref{sec:hsdpb}, we investigate those homogeneous 
projective bundles that are self-dual, i.e., equipped with
an isomorphism to their dual bundle; these correspond to 
principal bundles under the projective orthogonal or symplectic
groups. Here the main ingredients are the Heisenberg groups 
associated to symplectic vector spaces over the field with two
elements. Also, we introduce a geometric notion of indecomposability 
(which differs from the group-theoretic notion of 
$L$-indecomposability defined in \cite{BBN05}), and obtain 
a structure result for indecomposable homogeneous self-dual bundles 
(Proposition \ref{prop:indsd}). 

\medskip

Throughout this article, the base field $k$ is algebraically closed,
of arbitrary characteristic $p \geq 0$. Most of our results on 
$\bP^{n-1}$-bundles hold under the assumption that $n$ is not a 
multiple of $p$; indeed, the structure of commutative subgroup schemes 
of $\PGL_n$ is much more complicated when $p$ divides $n$ 
(see \cite{LMT09}). For the same reason, we only consider self-dual 
projective bundles in characteristic $\neq 2$. It would be interesting 
to extend our results to `bad' characteristics. 

\bigskip

\noindent
{\bf Acknowledgements}. This work originates in a series of lectures 
given at the Chennai Mathematical Institute in January 2011. I thank 
that institute and the Institute of Mathematical Sciences, Chennai, 
for their hospitaliy, and all the attendants of the lectures, especially
V. Balaji, P. Samuel and V. Uma, for their interest and stimulating 
questions. I also thank C. De Concini, P. Gille, F. Knop, and C. Procesi 
for very helpful discussions.

\bigskip

\noindent
{\bf Notation and conventions}.
We use the book \cite{DG70} by Demazure and Gabriel as a general
reference for group schemes. Our reference for abelian varieties 
is Mumford's book \cite{Mum70}; we generally follow its notation. 
In particular, the group law of an abelian variety $X$ is denoted 
additively and the multiplication by an integer $n$ is denoted by 
$n_X$, with kernel $X_n$. For any point $x \in X$, we denote by 
$T_x : X \to X$ the translation $y \mapsto x + y$. 
The dual abelian variety is denoted by $\wX$.

\section{Structure of homogeneous projective bundles}
\label{sec:str}

\subsection{Generalities on projective bundles}
\label{subsec:genpb}

Recall that a \emph{projective bundle} over a variety $X$ 
is a variety $P$ equipped with a proper flat morphism
\begin{equation}\label{eqn:proj}
f: P \longrightarrow X
\end{equation}
with fibers at all closed points isomorphic to projective space 
$\bP^{n-1}$ for some integer $n \geq 1$. Then $f$ is a 
$\bP^{n-1}$-bundle for the \'etale topology (see \cite[Sec.~I.8]{Gr68}).

Also, recall from [loc.~cit.] that the $\bP^{n-1}$-bundles are 
in a one-to-one correspondence with the torsors (or principal bundles)
\begin{equation}\label{eqn:prin}
\pi: Y \longrightarrow X
\end{equation}
under the projective linear group, $\PGL_n = \Aut(\bP^{n-1})$.
Specifically, $P$ is the associated bundle $Y \times^{\PGL_n} \bP^{n-1}$,
and $Y$ is the bundle of isomorphisms $X \times \bP^{n-1} \to P$
over $X$. Thus, any representation $\rho : \PGL_n \to \GL(V)$ defines 
the associated vector bundle $Y \times^{\PGL_n} V$ over $X$. 
The representation of $\PGL_n$ in the space $\M_n$ of $n \times n$ 
matrices by conjugation yields a `matrix bundle' on $X$; its 
sheaf of local sections is an Azumaya algebra of rank $n^2$ over $X$,
$$
\cA := (\pi_*(\cO_Y) \otimes \M_n)^{\PGL_n},
$$ 
viewed as a sheaf of non-commutative $\cO_X$-algebras over 
$\pi_*(\cO_Y)^{\PGL_n} = \cO_X$. In particular, $\cA$ defines 
a central simple algebra of degree $n$ over the function field 
$k(X)$. By \cite[Cor.~I.5.11]{Gr68}, 
the assignement $P \mapsto \cA$ yields a one-to-one correspondence
between $\bP^{n-1}$-bundles and Azumaya algebras of rank $n^2$. 
The quotient of $\cA$ by $\cO_X$ is the sheaf of local sections 
of the \emph{adjoint bundle} $\ad(P)$, the vector bundle associated with 
the adjoint representation of $\PGL_n$ in its Lie algebra 
$\mathfrak{pgl}_n$ (the quotient of the Lie algebra $\M_n$ by the scalar 
matrices). The correspondences between $\bP^{n-1}$-bundles, 
$\PGL_n$-torsors, and Azumaya algebras of rank $n^2$ preserve morphisms.
As a consequence, every morphism of $\bP^{n-1}$-bundles is an isomorphism.

There is a natural operation of \emph{product} on projective bundles:
to any $\bP^{n_i -1}$-bundles $f_i: P_i \to X$ ($i = 1,2$) with associated 
$\PGL_{n_i}$- bundles $\pi_i: Y_i \to X$, one associates the 
$\bP^{n_1 n_2 - 1}$-bundle 
$$
f : P_1 P_2 \longrightarrow X
$$ 
that corresponds to the $\PGL_{n_1 n_2}$-torsor
obtained from the $\PGL_{n_1} \times \PGL_{n_2}$-torsor
$$
\pi_1 \times \pi_2: Y_1 \times_X Y_2 \longrightarrow X
$$ 
by the extension of structure groups
$$
\PGL_{n_1} \times \PGL_{n_2} = \PGL(k^{n_1}) \times \PGL(k^{n_2})
\stackrel{\rho}{\longrightarrow} \PGL(k^{n_1} \otimes k^{n_2}) 
= \PGL_{n_1 n_2},
$$
where $\rho$ stems from the natural representation
$\GL(k^{n_1}) \times \GL(k^{n_2}) \to \GL(k^{n_1} \otimes k^{n_2})$.   
So $P_1P_2$ contains the fibered product $P_1 \times_X P_2$; 
it may be viewed as a global analogue of the Segre product of
projective spaces. The corresponding operation on Azumaya 
algebras is the tensor product (see \cite[Sec.~I.8]{Gr68}).

Likewise, any projective bundle $f : P \to X$ has a \emph{dual} bundle 
$$
f^* : P^* \longrightarrow X,
$$ 
where $P^*$ is the same variety as $P$, but the action of $\PGL_n$ is 
twisted by the automorphism arising from inverse transpose; then
$P^* = Y \times^{\PGL_n}(\bP^{n-1})^*$, where $(\bP^{n-1})^*$ denotes 
the dual projective space. The Azumaya algebra associated with $P^*$ 
is the opposite algebra $\cA^{\op}$. The assignement $P \mapsto P^*$ 
is of course contravariant, and the bi-dual $P^{**}$ comes with a
canonical isomorphism of bundles 
$P \stackrel{\cong}{\longrightarrow} P^{**}$.

Given a positive integer $n_1 \leq n$, a $\bP^{n_1 -1}$-\emph{sub-bundle} 
$f_1: P_1 \to X$ of the $\bP^{n-1}$-bundle (\ref{eqn:proj})
corresponds to a reduction of structure group of the associated 
$\PGL_n$-torsor (\ref{eqn:prin}) to a $\PGL_{n,n_1}$-torsor 
$\pi_1: Y_1 \to X$, where $\PGL_{n,n_1} \subset \PGL_n$ denotes 
the maximal parabolic subgroup that stabilizes a linear subspace 
$\bP^{n_1 - 1}$ of $\bP^{n - 1}$. Equivalently, the sub-bundle $P_1$ 
corresponds to a $\PGL_n$-equivariant morphism
$$
\gamma : Y \to \PGL_n / \PGL_{n,n_1} = \Gr_{n,n_1}
$$
(the Grassmannian parametrizing these subspaces). We have 
$P \cong Y_1 \times^{\PGL_{n,n_1}} \bP^{n - 1}$ and
$P_1 \cong Y_1 \times^{\PGL_{n,n_1}} \bP^{n_1 - 1}$
as bundles over $X$, where $\PGL_{n,n_1}$ acts on $\bP^{n_1 - 1}$
via its quotient $\PGL_{n_1}$.

Given two positive integers $n_1,n_2$ such that $n_1 + n_2 = n$,
a \emph{decomposition} of type $(n_1,n_2)$ of the $\bP^{n-1}$-bundle 
(\ref{eqn:proj}) consists of two $\bP^{n_i - 1}$-sub-bundles 
$f_i : P_i \to X$ ($i = 1,2$) which are disjoint (as subvarieties 
of $P$). This corresponds to a reduction of structure group of 
the $\PGL_n$-torsor (\ref{eqn:prin}) to a torsor 
$\pi_{12}: Y_{12} \to X$ under the maximal Levi subgroup 
$$
\P(\GL_{n_1}\times\GL_{n_2}) = \PGL_{n,n_1} \cap \PGL_{n,n_2} 
\subset \PGL_n
$$  
that stabilizes two disjoint linear subspaces $\bP^{n_i - 1}$ 
of $\bP^{n-1}$ ($i = 1,2$). Then 
$$
P_i = Y_{12} \times^{\P(\GL_{n_1}\times\GL_{n_2})} \bP^{n_i -1}
$$
for $i = 1,2$, where $\P(\GL_{n_1}\times\GL_{n_2})$ acts on each
$\bP^{n_i -1}$ via its quotient $\PGL_{n_i}$. The decompositions of type 
$(n_1,n_2)$ correspond to the $\PGL_n$-equivariant morphisms
\begin{equation}\label{eqn:dec}
\delta : Y \longrightarrow \PGL_n/\P(\GL_{n_1} \times \GL_{n_2})
\end{equation}
to the variety of decompositions.

If the bundle (\ref{eqn:proj}) admits no decomposition, 
then we say of course that it it \emph{indecomposable}. 
Equivalently, the associated torsor (\ref{eqn:prin}) admits no 
reduction of structure group to a proper Levi subgroup. 

When $P$ is the projectivization $\bP(E)$ of a vector bundle 
$E$ over $X$, the sub-bundles of $P$ correspond bijectively to those 
of $E$, and the decompositions of $P$, to the splittings 
$E = E_1 \oplus E_2$ of vector bundles. Also, note that 
$\bP(E) \, \bP(F) = \bP(E \otimes F)$ and $\bP(E)^* = \bP(E^*)$ 
with an obvious notation.

\subsection{Homogeneous projective bundles}
\label{subsec:hpb}

From now on, $X$ denotes a fixed abelian variety, $f: P \to X$ 
a $\bP^{n-1}$-bundle, and $\pi : Y \to X$ the corresponding 
$\PGL_n$-torsor. Then $P$ is a nonsingular projective variety
and $f$ is its Albanese morphism. In particular, $f$ is uniquely 
determined by the variety $P$.

Since $P$ is complete, its automorphism functor is represented 
by a group scheme $\Aut(P)$, locally of finite type. Moreover, 
we have a homomorphism of group schemes
$$
f_* : \Aut(P) \longrightarrow \Aut(X)
$$
with kernel the subgroup scheme $\Aut_X(P) \cong \Aut_X^{\PGL_n}(Y)$ 
of bundle automorphisms. Also, $\Aut_X(P)$ is affine of finite type,
and its Lie algebra is $H^0(X, \ad(P))$ (see e.g. \cite[Sec.~4]{Br11}
for these results). 

We say that a $\bP^{n-1}$-bundle (\ref{eqn:proj}) is 
\emph{homogeneous}, if the image of $f_*$ contains the subgroup
$X \subset \Aut(X)$ of translations; equivalently, the bundle $P$
is isomorphic to its pull-backs under all translations.
This amounts to the vector bundle $\ad(P)$ being homogeneous
(see \cite[Cor.~2.15]{Br12}; if $P$ is the projectivization of 
a vector bundle, this follows alternatively from \cite[Thm.~5.8]{Muk78}).
 
The structure of homogeneous projective bundles is described by the 
following:

\begin{theorem}\label{thm:str}
{\rm (i)} A $\bP^{n-1}$-bundle $f: P \to X$ is homogeneous 
if and only if there exist an exact sequence of group schemes 
\begin{equation}\label{eqn:ext}
\CD
1 @>>> H @>>> G @>{\gamma}>> X @>>> 1,
\endCD
\end{equation}
where $G$ is anti-affine (i.e., $\cO(G) = k$), and a faithful 
homomorphism $\rho : H \longrightarrow \PGL_n$ such that $P$ is 
the associated bundle $G \times^H \bP^{n-1} \to G/H = X$, 
where $H$ acts on $\bP^{n-1}$ via~$\rho$. 

Then the exact sequence (\ref{eqn:ext}) is unique; the group scheme
$G$ is smooth, connected, and commutative (in particular, $H$ is
commutative), and the projective representation $\rho$ is unique 
up to conjugacy in $\PGL_n$.  
Moreover, the corresponding $\PGL_n$-torsor
is the associated bundle $G \times^H \PGL_n \to X$, and the 
corresponding Azumaya algebra satisfies
\begin{equation}\label{eqn:azu}
\cA \cong (\gamma_*(\cO_G) \otimes \M_n)^H
\end{equation}
as a sheaf of algebras over $\gamma_*(\cO_G)^H \cong \cO_X$.

\smallskip

\noindent
{\rm (ii)} For $P$ as in {\rm (i)}, we have an isomorphism
\begin{equation}\label{eqn:aut}
\Aut_X(P) \cong \PGL_n^H
\end{equation}
(the centralizer of $H$ in $\PGL_n$). As a consequence, 
\begin{equation}\label{eqn:lie}
H^0(X, \ad(P)) = \fpgl_n^H.
\end{equation}

\smallskip

\noindent
{\rm (iii)} The homogeneous projective sub-bundles of $P$ are 
exactly the bundles $G \times^H S \to X$, where 
$S \subset \bP^{n-1}$ is an $H$-stable linear subspace.

\smallskip

\noindent
{\rm (iv)} Any decomposition of $P$ consists of homogeneous
sub-bundles.
\end{theorem} 

\begin{proof}
(i) follows readily from Theorem 3.1 in \cite{Br12}, and (ii)
from Proposition 3.6 there.

(iii) Let $f_1: P_1 \to X$ be a projective sub-bundle, and consider
the corresponding reduction of structure group of the 
$\PGL_n$-torsor $Y$ to a $\PGL_{n,n_1}$-torsor 
$\pi_1: Y_1 \to X$. If $f_1$ is homogeneous, then by 
\cite[Thm.~3.1]{Br12} again, we have a $\PGL_{n,n_1}$-equivariant 
isomorphism
$$
Y_1 \cong G_1 \times^{H_1} \PGL_{n,n_1}
$$
for some exact sequence $0 \to H_1 \to G_1 \to X \to 0$ with $G_1$
anti-affine, and some faithful homomorphism 
$\rho_1 : H_1 \to \PGL_{n,n_1}$. Thus,
$$
Y \cong Y_1 \times^{\PGL_{n,n_1}} \PGL_n \cong G_1 \times^{H_1} \PGL_n
$$
equivariantly for the action of $\PGL_n$. By the uniqueness in (i),
it follows that $G_1 = G$ and $H_1 = H$; hence $P_1 = G \times^H S$
for some $H$-stable linear subspace $S \subset \bP^{n-1}$.

Conversely, any $H$-stable linear subspace obviously yields 
a homogeneous projective sub-bundle.

(iv) A decomposition of $P$ of type $(n_1,n_2)$ corresponds to
a $\PGL_n$-equivariant morphism 
$\delta : Y \to \PGL_n/\P(\GL_{n_1} \times \GL_{n_2})$. Since the 
variety $\PGL_n/\P(\GL_{n_1} \times \GL_{n_2})$ is affine, 
the corresponding reduction of structure group 
$\pi_{12}: Y_{12} \to X$ is homogeneous by \cite[Prop.~2.8]{Br12}. 
Thus, the associated bundles $P_1, P_2$ are homogeneous as well.
\end{proof}

\begin{remark}\label{rem:prod}
Let $P_i$ ($i = 1, 2$) be homogeneous bundles corresponding 
to extensions $1 \to H_i \to G_i \to X \to 1$ and projective 
representations $\rho_i : H_i \to \PGL_{n_i}$. Then the 
$\PGL_{n_1n_2}$-torsor that corresponds to $P_1 P_2$ is the 
associated bundle
$$
(G_1 \times_X G_2) \times^{H_1 \times H_2} \PGL_{n_1 n_2} 
\longrightarrow (G_1 \times_X G_2)/(H_1 \times H_2) = X,
$$ 
where the homomorphism $H_1 \times H_2 \to \PGL_{n_1n_2}$
is given by the tensor product $\rho_1 \otimes \rho_2$.
Thus, $P_1 P_2$ is the homogeneous bundle classified by the
extension $1 \to H \to G \to X \to 1$, where 
$G \subset G_1 \times_X G_2$ denotes the largest anti-affine subgroup
and $H = (H_1 \times H_2) \cap G$, and by the projective 
representation $(\rho_1 \otimes \rho_2)\vert_H$. 

As a consequence, the $m$th power $P^m$ corresponds to the same 
extension as $P$ and to the $m$th tensor power of its projective 
representation. Likewise, the dual of a homogeneous bundle is the 
homogeneous bundle associated with the same extension and with 
the dual projective representation.
\end{remark}

The anti-affine algebraic groups are classified in \cite{Br09} 
and independently \cite{SS09}, and the anti-affine extensions 
(\ref{eqn:ext}) in \cite[Sec.~3.3]{Br12}. We now describe 
the other ingredients of Theorem \ref{thm:str}, i.e., 
the commutative subgroup schemes $H \subset \PGL_n$ up to conjugacy. 
Every such subgroup scheme has a unique decomposition
$$
H = H_u \times H_s,
$$
where $H_u$ is unipotent and $H_s$ is diagonalizable.
Thus, $H_s$ sits in an exact sequence
$$
1 \longrightarrow H_s^0 \longrightarrow H_s \longrightarrow F 
\longrightarrow 1,
$$ 
where $H_s^0$ is a connected diagonalizable group scheme
(the neutral component of $H_s$), 
and the group of components $F$ is finite, diagonalizable and 
of order prime to $p$ (in particular, $F$ is smooth); 
this exact sequence is unique and splits non-canonically. 
In turn, $H_s^0$ is an extension of a finite diagonalizable group 
scheme of order a power of $p$, by a torus (the reduced neutral 
component); this extension is also unique and splits non-canonically.

Denote by $\tH \subset \GL_n$ the preimage of $H \subset \PGL_n$.
This yields a central extension
\begin{equation}\label{eqn:cext}
1 \longrightarrow \bG_m \longrightarrow \tH \longrightarrow H 
\longrightarrow 1,
\end{equation}
where the multiplicative group $\bG_m$ is viewed as the group 
of invertible scalar matrices. We say that $\tH$ is the 
\emph{theta group} of $H$, and define similarly $\tHu,\tHs$
and $\tHso$ (the latter is the neutral component of $\tHs$).

Given two $S$-valued points $\tx, \ty$ of $\tH$, 
where $S$ denotes an arbitrary scheme, the commutator 
$\tx \ty \tx^{-1} \ty^{-1}$ is a $S$-valued point of 
$\bG_m$ and depends only on the images of $\tx, \ty$ in $H$. 
This defines a morphism
\begin{equation}\label{eqn:comm}
e: H \times H \longrightarrow \bG_m
\end{equation}
which is readily seen to be bilinear (i.e., we have 
$e(xy,z) = e(x,z) \, e(y,z)$ 
and $e(x,yz) = e(x,z) \, e(y,z)$ for all $S$-valued points
$x,y,z$ of $H$) and alternating (i.e., $e(x,x) = 1$ for all $x$). 
We say that $e$ is the \emph{commutator pairing} of the extension
(\ref{eqn:cext}). 

Note that the dual bundle $P^*$ has pairing $e^{-1}$; moreover, 
the power $P^m$, where $m$ is a positive integer, has pairing $e^m$. 
 
The center $Z(\tH)$ sits in an exact sequence of group schemes
\begin{equation}\label{eqn:zext}
1 \longrightarrow \bG_m \longrightarrow Z(\tH) \longrightarrow H^{\perp}
\longrightarrow 1,
\end{equation}
where the $S$-valued points of $H^{\perp}$ are those points of 
$H$ such that $e(x,y) = 1$ for all $S'$-valued points $y$ of 
$H$ and all schemes $S'$ over $S$. In particular, 
\emph{$\tH$ is commutative if and only if $e = 1$}.

We now show that the obstruction for being the projectivization 
of a homogeneous vector bundle is just the commutator pairing.
The obstruction for being the projectivization of an arbitrary 
vector bundle will be determined in Theorem \ref{thm:obs}.

\begin{proposition}\label{prop:homvec}
With the above notation, the following conditions are equivalent:

\smallskip

\noindent
{\rm (i)} $P$ is the projectivization of a homogeneous vector bundle.

\smallskip

\noindent
{\rm (ii)} The extension (\ref{eqn:cext}) splits.

\smallskip

\noindent
{\rm (iii)} $e = 1$.
\end{proposition}

\begin{proof}
(i)$\Rightarrow$(ii) By \cite[Thm.~3.1]{Br12}, 
any homogeneous vector bundle $E$ of rank $n$ over $X$ is of the form
$G \times^H k^n \longrightarrow G/H = X$
for some anti-affine extension $1 \to H \to G \to X \to 1$ and some
faithful representation $\sigma : H \to \GL_n$. Since $H$ is commutative,
$k^n$ contains eigenvectors of $H$; thus, twisting $\sigma$ by a character 
of $H$ (which does not change the projectivization $\bP(E)$), we may assume 
that $k^n$ contains non-zero fixed points of $H$. Then $\sigma$ defines a 
faithful projective representation $\rho : H \to \PGL_n$. Hence $G$ and 
$\rho$ are the data associated with the homogeneous projective bundle 
$\bP(E) \to X$, and $\sigma$ splits the extension (\ref{eqn:cext}).

(ii)$\Rightarrow$(i) Any splitting of that extension yields 
a homomorphism $\sigma : H \to \GL_n$ that lifts $\rho$. 
Then the associated bundle $G \times^H k^n \to X$ is a homogeneous 
vector bundle with projectivization $P$.

(ii)$\Rightarrow$(iii) is obvious. Conversely, if $e = 1$, 
then $\tH$ is commutative. It follows that
$\tH \cong U \times \tHs$, where the unipotent part $U$ is
isomorphic to $H_u$ via the homomorphism $\tH \to H$, and
$\tHs$ sits in an exact sequence of diagonalizable group schemes
$1 \to \bG_m \to \tHs \to H_s \to 1$.
But every such sequence splits, since so does the dual exact sequence
of character groups.
\end{proof}

Next, we obtain a very useful structure result for $H$ under 
the assumption that $n$ is not divisible by the characteristic:

\begin{proposition}\label{prop:theta}
Keep the above notation, and assume that $(n,p) =1$.

\smallskip

\noindent
{\rm (i)} The extension $1 \to \bG_m \to \tHu \to H_u \to 1$ 
has a unique splitting, and the corresponding lift of $H_u$
(that we still denote by $H_u$) is central in $\tH$. Also, 
the extension $1 \to \bG_m \to \tHso \to H_s^0 \to 1$ 
splits non-canonically and $\tHso$ is central in $\tH$. 

\smallskip

\noindent
{\rm (ii)} We have canonical decompositions of group schemes 
$$
\tH = H_u \times \tHs, \quad Z(\tH) = H_u \times Z(\tHs).
$$
Moreover, $Z(\tHs)$ is diagonalizable and sits in an exact sequence
$$
1 \longrightarrow \tHso \longrightarrow Z(\tHs) \longrightarrow
F^{\perp} \longrightarrow 1
$$
which splits non-canonically.

\smallskip

\noindent
{\rm (iii)} The commutator pairing $e$ factors through a 
bilinear alternating morphism
\begin{equation}\label{eqn:cpf} 
e_F : F \times F \longrightarrow \bG_m.
\end{equation}
\end{proposition}

\begin{proof}
Since any commutator has determinant $1$, 
we see that $e$ takes values in the subgroup scheme 
$\mu_n = \bG_m \cap \SL_n$ of $n$th roots of unity. In other 
terms, $e$ factors through the pairing
$$
se : H \times H \longrightarrow \mu_n
$$
defined by the central extension
$$ 
1 \longrightarrow \mu_n \longrightarrow S\tH 
\longrightarrow H \longrightarrow 1,
$$
where $S\tH := H \cap \SL_n$. Note that $\mu_n$ is smooth
by our assumption on $n$. Moreover, $se$ restricts trivially to
$nH \times H$, where $nH$ denotes the image of the 
multiplication by $n$ in the commutative group scheme $H$. 

We claim that $H_u \subset nH$. This is clear if $p = 0$,
since $H_u$ is then isomorphic to the additive group of
a vector space. If $p \geq 1$, then the commutative
unipotent group scheme $H_u$ is killed by some power of $p$.
Using again the assumption that $(n,p) = 1$, it follows
that $H_u = n H_u \subset n H$.

By that claim, $se$ restricts trivially to $H_u \times H$,
and hence $\tHu \subset Z(\tH)$; in particular, $\tHu$ is
commutative. Thus, $\tHu \cong H_u \times \bG_m$; this
proves the assertion about $H_u$.

We already saw that the extension 
$1 \to \bG_m \to \tHso \to H_s^0 \to 1$ splits. 
Also, $H_s^0 \cong T \times E$, where $T$ is a torus 
(the reduced neutral component), and $E$ is a finite group scheme
killed by some power of $p$. As above, it follows that 
$H_s^0 \subset n H$, and that $\tHso$ is central in $\tH$. 
This completes the proof of (i). 

The decompositions in (ii) are direct consequences of (i).
The assertion on $Z(\tHs)$ follows from the exact sequence 
$1  \to \tHso \to \tHs \to F \to 1$, since 
$\tHso \subset Z(\tHs)$. 
Finally, (iii) also follows readily from (i).
\end{proof}

\begin{remark}\label{rem:form}
With the notation and assumptions of the above proposition,
the group scheme $\Aut_X(P)$ is smooth, as follows from 
the isomorphism (\ref{eqn:aut}) together with 
\cite[Thm.~1.1]{He10}. Moreover, $\Aut(P)$ is smooth as well: 
indeed, we have an exact sequence of group schemes
$$
\CD
1 @>>> \Aut_X(P) @>>> \Aut(P) @>{f_*}>> \Aut_P(X) @>>> 1,
\endCD
$$
where $\Aut_P(X)$ is a subgroup scheme of $\Aut(X)$ containing the
group $X$ of translations. Since $\Aut(X) = X \ltimes \Aut_{\gp}(X)$,
where the group scheme of automorphisms of algebraic groups 
$\Aut_{\gp}(X)$ is \'etale (possibly infinite), it follows that 
$\Aut_P(X)$ is smooth, and hence so is $\Aut(X)$.
\end{remark}

\subsection{Non-degenerate theta groups}
\label{subsec:ndtg}

As in the above subsection, we consider a commutative subgroup
scheme $H \subset \PGL_n$ and the associated theta group 
$\tH \subset \GL_n$; we assume that $(n,p) = 1$.
 
We say that $\tH$ is \emph{non-degenerate} if $Z(\tH) = \bG_m$. 
By Proposition \ref{prop:theta}, this is equivalent to the assertions
that $H$ is a finite commutative group of order prime to $p$, 
and the homomorphism 
\begin{equation}\label{eqn:dual}
\epsilon : H \longrightarrow \cX(H), \quad x \longmapsto
(y \mapsto e(x,y))
\end{equation} 
is faithful, where $\cX(H) := \Hom_{\gp}(H,\bG_m)$ 
denotes the character group of $H$. It follows that $\epsilon$ 
is an isomorphism.
 
We now recall from \cite[Sec.~1]{Mum66} the structure of 
non-degenerate theta groups. Choose a subgroup $K \subset H$ 
that is totally isotropic for the commutator pairing $e$, 
and maximal with this property. Then 
$$
\tH \cong \bG_m \times K \times \cX(K),
$$ 
where the  group law on the right-hand side is given by
\begin{equation}\label{eqn:law}
(t,x,\chi) \cdot (t',x',\chi') = (tt' \chi'(x), x + x', \chi+ \chi'),
\end{equation}
the group laws on $K$ and $\cX(K)$ being denoted additively.
Such a group is called the \emph{Heisenberg group} associated with 
the finite group $K$; we denote it by $\cH(K)$ and identify the group 
$K$ (resp. $\cX(K)$) with its lift $\{ 1 \} \times K \times \{ 0 \}$ 
(resp. $\{ 1 \} \times \{ 0 \} \times \cX(K)$) in $\cH(K)$.

Also, recall that $\cH(K)$ has a unique irreducible representation 
on which $\bG_m$ acts via $t \mapsto t \, \id$: the 
\emph{standard representation} (also called the
\emph{Schr\"odinger representation}) 
in the space $\cO(K)$ of functions on $K$ with values in $k$, 
on which $\tH$ acts via
$$
((t,x,\chi) \cdot f)(y) := t \, \chi(y) \, f(x + y).
$$
The corresponding commutator pairing $e$ is given by
$$
e((x,\chi), (x',\chi')) := \chi'(x) \, \chi(x')^{-1}.
$$
In particular, the standard representation $W(K)$ contains a unique line of 
$K$-fixed points and has dimension $n= \#(K)$; moreover, the group $H$ 
is killed by $n$ and has order $n^2$. 
Any finite-dimensional representation $V$ of $\cH(K)$ on which $\bG_m$ acts 
by scalar multiplication is a direct sum of $m$ copies of $W(K)$, where 
$m := \dim(V^K)$. Such a representation is called \emph{of weight $1$.}

For later use, we record the following result, which is well-known 
in the setting of theta structures on ample line bundles over complex
abelian varieties (see \cite[Lem.~6.6.6 and Exer.~6.10.14]{BL04}):

\begin{lemma}\label{lem:cent}
Assume that $(n,p) = 1$ and let $\tH \subset \GL_n$ be a 
non-degenerate theta group.

\smallskip

\noindent
{\rm (i)} The algebra $\M_n$ has a basis $(u_h)_{h \in H}$ such that
every $u_h$ is an eigenvector of $H$ (acting by conjugation) with
weight $\epsilon(h)$, and
$$
u_{x,\chi} \, u_{x',\chi'} = \chi'(x) \,u_{x + x', \chi + \chi'}
$$
for all $h = (x, \chi)$ and $h' = (x',\chi')$ in $H = K \times \cX(K)$. 
In particular, the representation of $H$ in $\M_n$ by conjugation is
isomorphic to the regular representation.

\smallskip

\noindent
{\rm (ii)} The centralisers of $\tH$ in $\GL_n$ and of $H$ in 
$\PGL_n$ satisfy 
$$
\GL_n^{\tH} = \bG_m, \quad \PGL_n^H = H.
$$
Moreover, the normalizers sit in exact sequences 
$$
1 \longrightarrow \bG_m  \longrightarrow N_{\GL_n}(\tH)
\longrightarrow N_{\PGL_n}(H) \longrightarrow 1,
$$
$$
1 \longrightarrow H \longrightarrow N_{\PGL_n}(H)
\longrightarrow \Aut(H,e) \longrightarrow 1.
$$
Also, we have an isomorphism
\begin{equation}\label{eqn:norm}
\Aut^{\bG_m}(\tH) \cong N_{\PGL_n}(H).
\end{equation}
\end{lemma}

\begin{proof}
(i) We may view $H$ as a subset of $\M_n$ via 
$(x,\chi) \mapsto u_{x,\chi} := (1,x,\chi) \in \tH \subset \GL_n$. 
Then the assertions follow readily from the formula (\ref{eqn:law})
for the group law of $\tH$. 

(ii) By Schur's lemma, we have $\GL_n^{\tH} = \bG_m$; this yields
the first exact sequence. 

In view of (i), the fixed points of $H$ acting on $\bP(M_n)$ by
conjugation are exactly the points of $H \subset \PGL_n$; thus,
$\PGL_n^H = H$. To obtain the second exact sequence, it suffices 
to show that the image in $\Aut(H)$ of $N_{\PGL_n}(H)$ equals
$\Aut(H,e)$. But if $g \in \PGL_n$ normalizes $H$, then one readily 
checks that the conjugation $\Int(g)\vert_H$ preserves the pairing 
$e$. Conversely, let $g \in \Aut(H,e)$; then composing the inclusion 
$\rho: H \to \PGL_n$ with $g$, we obtain a projective representation 
$\rho_g$ with the same commutator pairing. Thus, $\rho_g$ lifts to 
a representation $\trho_g : \tH \to \GL_n$ which is isomorphic to 
the standard representation. It follows that $g$ extends to the 
conjugation by some $\tg \in \GL_n$ that normalizes $H$.

The isomorphism (\ref{eqn:norm}) follows similarly from the 
fact that the standard representation is the unique irreducible
representation of weight $1$.
\end{proof}

Returning to an arbitrary theta group $\tH \subset \GL_n$, 
we now describe the representation of $\tH$ in $k^n =: V$. 
Consider the decomposition
\begin{equation}\label{eqn:eigen}
V = \bigoplus_{\lambda} V_{\lambda}
\end{equation}
into weight spaces of the diagonalizable group $Z(\tHs)$, where $\lambda$
runs over the characters of weight $1$ of that group (those that restrict 
to the identity character of $\bG_m$). By Proposition \ref{prop:theta}, 
each $V_{\lambda}$ is stable under $\tH$.

\begin{proposition}\label{prop:rep}
With the above notation, each quotient $\tHs/\ker(\lambda)$ 
is isomorphic to the Heisenberg group $\cH(K/F^{\perp})$, 
where $K$ denotes a maximal totally isotropic subgroup scheme
of $F$ relative to $e_F$.

Moreover, we have an isomorphism of representations of 
$\tH \cong H_u \times \tHs$:
$$
V_{\lambda} \cong U_{\lambda} \otimes W(K/F^{\perp}),
$$
where $U_{\lambda}$ is a representation of $H_u$ and 
$W(K/F^{\perp})$ is the standard representation of 
$\tHs/\ker(\lambda)$.
\end{proposition}

\begin{proof}
Note that $\lambda$ yields a splitting of (\ref{eqn:zext}),
and an isomorphism $Z(\tHs)/\ker(\lambda) \cong \bG_m$. Also, 
$\tHs/Z(\tHs) \cong \tH/Z(\tH) \cong F/F^{\perp}$ 
by Proposition \ref{prop:theta}. Thus, the exact sequence
$$
1 \longrightarrow Z(\tHs)/\ker(\lambda) \longrightarrow 
\tHs/\ker(\lambda) \longrightarrow \tHs/Z(\tHs) 
\longrightarrow 1
$$
may be identified with the central extension 
$$ 
1 \longrightarrow \bG_m \longrightarrow \tHs/\ker(\lambda)
\longrightarrow F/F^{\perp} \longrightarrow 1
$$ 
and the corresponding commutator pairing is induced by $e_F$.
This shows that $\tHs/\ker(\lambda)$ is a non-degenerate theta group.
Now the first assertion follows from the structure of these groups. 

Also, $V_{\lambda}$ is a representation of $\tHs/\ker(\lambda)$ 
on which the center $\bG_m$ acts with weight $1$,
and hence a direct sum of copies of the standard representation.
This implies the second assertion in view of Proposition 
\ref{prop:theta} again.
\end{proof}

\begin{corollary}\label{cor:deg}
With the above notation, the representation of $\tH$ in $V$ is an 
iterated extension of irreducible representations of the same dimension,
\begin{equation}\label{eqn:hindex}
d := [K:F^{\perp}] = \sqrt{[F:F^{\perp}]}= \sqrt{[H:H^{\perp}]}.
\end{equation}
In particular, $n$ is a multiple of $d$, with equality if and only if 
$\tH$ is a Heisenberg group acting via its standard representation.
\end{corollary}
 
We say that $d$ is the \emph{homogeneous index} of the bundle 
(\ref{eqn:proj}); this is the minimal rank of a homogeneous sub-bundle 
of $P$ in view of Theorem \ref{thm:str}. (One can show that the 
homogeneous index of $P$ is a multiple of the index of the associated
central simple algebra over $k(X)$). 
Note that $F/F^{\perp}$ is killed by $d$, and hence $e_F^d = 1$. 
In view of Proposition \ref{prop:homvec}, it follows that 
the $d$th power $P^d$ 
\emph{is the projectivization of a homogeneous vector bundle}.

\begin{proposition}\label{prop:indec}
With the notation and assumptions of Proposition \ref{prop:theta}, 
the following assertions are equivalent for a homogeneous 
$\bP^{n-1}$-bundle $f: P \to X$:

\smallskip

\noindent
{\rm (i)} $P$ is indecomposable.

\smallskip

\noindent
{\rm (ii)} The associated representation $\trho: \tH \to \GL_n$
is indecomposable.

\smallskip

\noindent
{\rm (iii)} $\tHs$ is a Heisenberg group and $V \cong U \otimes W$ 
as representations of $H \cong H_u \times \tHs$, where $U$ is 
an indecomposable representation of $H_u$ and $W$ is the standard 
representation of $\tHs$.

\smallskip

\noindent
{\rm (iv)} The neutral component $\Aut^0_X(P)$ is unipotent.
\end{proposition}

\begin{proof}
(i)$\Rightarrow$(ii) is obvious, and (ii)$\Rightarrow$(i) 
follows from Theorem \ref{thm:str} (iv). 

(ii)$\Leftrightarrow$(iii) is a direct consequence of Proposition 
\ref{prop:rep}.

(iii)$\Rightarrow$(iv) Since $(n,p) = 1$, we have 
$\M_n = k \id \oplus \fpgl_n$ as representations of $\PGL_n$ acting
by conjugation. In view of (\ref{eqn:lie}), this yields
$$
\Lie \, \Aut_X(P) = \M_n^H/k \, \id = \End^H(U \otimes W)/k \, \id.
$$
Moreover, $\End^H(U \otimes W) \cong \End^{H_u}(U)$ 
by Schur's lemma, and hence 
$$
\Lie \, \Aut_X(P) \cong \End^{H_u}(U)/k \, \id.
$$
This isomorphism of Lie algebras arises from the natural homomorphism
$$
\GL(U)^{H_u}/\bG_m \, \id \longrightarrow \Aut_X(P).
$$ 
Since $\Aut_X(P)$ is smooth (Remark \ref{rem:form}), we see that
its neutral component is a quotient of $\GL(U)^{H_u}/\bG_m \, \id$.
But the latter group is unipotent, since $U$ is indecomposable.

(iv)$\Rightarrow$(iii) Observe that the weight space decomposition 
(\ref{eqn:eigen}) is trivial: otherwise, $\Aut_X(P)$ contains a 
copy of $\bG_m$ that fixes some weight space pointwise and acts 
by scalar multiplication on all the other weight spaces. 
Thus, $V \cong U \otimes W$, where $W$ is irreducible. Moreover, 
$U$ is indecomposable; otherwise, $\Aut_X(P)$ contains a copy of 
$\bG_m$ by the above argument.   
\end{proof}

\begin{remarks}\label{rem:ks}
(i) The results of this subsection do not extend readily to the case 
where $p$ divides $n$: for instance, there exists a non-degenerate theta 
group $\tH \subset \GL_p$ with $H$ unipotent and local.
Consider indeed the group scheme $\alpha_p$ (the kernel of the
$p$th power map of $\bG_a$) and the duality pairing
$$
u : \alpha_p \times \alpha_p \longrightarrow \bG_m, \quad
(x,y) \longmapsto \sum_{i=0}^{p-1} \frac{x^i}{i!}.
$$
This yields a bilinear alternating pairing $e$ on 
$H := \alpha_p \times \alpha_p$ via
$$
e((x,y),(x',y')) := u(x, y') \, u(x',y)^{-1}.
$$
Then we may take for $\tH$ the associated Heisenberg group scheme
(with $K = \alpha_p \times \{ 0 \}$ and $\cX(K) = \{ 0 \} \times \alpha_p$),
equipped with its standard representation in $\cO(\alpha_p) \cong k^p$.

Note that the above group scheme $H$ is contained in an abelian surface,
the product of two supersingular elliptic curves. More generally,
any finite commutative group scheme is contained in some abelian
variety (see \cite[Sec.~15.4]{Oo66}).

\smallskip

\noindent
(ii) For an arbitrary homogeneous projective bundle $P$, each representation 
$U_{\lambda}$ (with the notation of Proposition \ref{prop:rep}) 
is a direct sum of indecomposable representations with multiplicities;
moreover, these indecomposable summands and their multiplicities
are uniquely determined up to reordering, in view of the Krull-Schmidt 
theorem. Thus, the representation of $\tH$ in $V$ decomposes into 
a direct sum (with multiplicities) of tensor products $U \otimes W$, 
where $U$ is an indecomposable representation of $H_u$ and $W$ is
an irreducible representation of $\tHs$. 

Let $L \subset \PGL_n$ denote the stabilizer of such a decomposition. 
Then $L$ is a Levi subgroup, uniquely determined up to conjugation; 
moreover, the $\PGL_n$-torsor $\pi : Y \to X$  
admits a reduction of structure group to an $L$-torsor
$\pi_L : Y_L \to X$. Arguing as in the proof of (iii)$\Rightarrow$(iv) 
above, one may check that the natural homorphism
$Z(L) \to \Aut^L_X(Y_L)$ (where $Z(L)$ denotes the center of $L$,
and $\Aut^L_X(Y_L)$ the group of bundle automorphisms of $Y_L$)
yields an isomorphism of the reduced neutral component $Z(L)^0_{\red}$ 
to a maximal torus of $\Aut^L_X(Y_L)$. Thus, the torsor
$\pi_L : Y_L \to X$ is $L$-\emph{indecomposable} in the sense of 
\cite[Def.~2.1]{BBN05}. Moreover, this torsor is the unique 
reduction of $\pi : P \to X$ to an $L$-indecomposable torsor for a Levi 
subgroup, by [loc.~cit., Thm.~3.4] (the latter result is obtained
there in characteristic zero, and generalized to arbitrary characteristics
in \cite{BBNP06}; see also \cite{BBN06}).

Conversely, the equivalence of the above statements (i) and (iv) 
follows from the results of \cite{BBN05, BBNP06} in view of 
the smoothness of $\Aut_X(P)$.
\end{remarks}

\section{Irreducible bundles}
\label{sec:irr}

Throughout this section, we consider $\bP^{n-1}$-bundles 
$f : P \to X$, and call them \emph{bundles} for simplicity;
we still assume that $(n,p) = 1$.

\subsection{Structure and characterizations}
\label{subsec:st}

We say that a homogeneous bundle $P$ is \emph{irreducible} 
if so is the projective representation $\rho : H \to \PGL_n$ 
associated with $P$ via Theorem \ref{thm:str}. 
By Proposition \ref{prop:rep}, this means that the theta group
$\tH$ is a Heisenberg group acting on $k^n$ via its standard
representation.

We now parametrize the irreducible homogeneous bundles, and describe
the corresponding Azumaya algebras as well as the adjoint bundles and 
automorphism groups:

\begin{proposition}\label{prop:irraut}
{\rm (i)} The irreducible homogeneous $\bP^{n-1}$-bundles are classified 
by the pairs $(H,e)$, where $H \subset X_n$ is a subgroup of 
order $n^2$ and $e : H \times H \to \bG_m$ is a non-degenerate 
alternating pairing. In particular, such bundles exist for any given
$n$, and they form only finitely many isomorphism classes.

\smallskip

\noindent
{\rm (ii)} For the bundle $P$ corresponding to $(H,e)$, the associated
Azumaya algebra $\cA$ admits a grading by the group $H$, namely,
$$
\cA \cong \bigoplus_{\cL \in H} \cL,
$$
where each element of $H \subset \wX$ is viewed as an invertible
sheaf on $X$. In particular, we have a decomposition
$$
\ad(P) \cong \bigoplus_{\cL \in H, \cL \neq 0} \cL.
$$
 
\smallskip

\noindent
{\rm (iii)} For $P$ as in {\rm (ii)}, we have $\Aut_X(P) \cong H$.
Moreover, 
the neutral component $\Aut^0(P)$ is the extension of $X$ by $H$, dual 
to the inclusion $\cX(H) \cong H \subset \wX$, and $\Aut(P)/\Aut^0(P)$ 
is isomorphic to the subgroup of $\Aut_{\gp}(X) \cong \Aut_{\gp}(\wX)$ 
that preserves $H$ and $e$.
\end{proposition}

\begin{proof}
(i) By the results of Section \ref{sec:str}, the irreducible homogeneous 
bundles are classified by the pairs consisting of an isogeny
$1 \to H \to G \to X \to 1$ and a non-degenerate alternating pairing
$e$ on $H$; then $e$ provides an isomorphism of $H$ with its character 
group. The assertion now follows from duality of isogenies.

(ii) follows from the isomorphism of $\cO_X$-algebras 
(\ref{eqn:azu}) together with the isomorphism of $\cO_X$-$H$-algebras 
$\gamma_*(\cO_G) \cong \bigoplus_{\cL \in \cX(H)} \cL$ and with 
the decomposition $\M_n \cong \bigoplus_{h \in H} k \, u_h$ 
obtained in Lemma \ref{lem:cent} (iii).

(iii) Combining the isomorphism (\ref{eqn:aut}) and Lemma 
\ref{lem:cent} (i), we see that the natural map $H \to \Aut_X(P)$ 
is an isomorphism. In view of the commutative diagram with exact rows
$$
\CD
1 @>>> H @>>> G @>>> X @>>> 1 \\
& & @VVV @VVV @VVV \\
1 @>>> \Aut_X(P) @>>> \Aut(P) @>{f_*}>> \Aut(X) \\
\endCD
$$
and of the isomorphism $\Aut(X) \cong X \ltimes \Aut_{\gp}(X)$,
where $\Aut_{\gp}(X)$ is \'etale, it follows that the natural 
map $G \to \Aut^0(P)$ is an isomorphism as well. The structure of
$\Aut(P)/\Aut^0(P)$ follows from Theorem \ref{thm:str} together with 
Lemma \ref{lem:cent} (ii). 
\end{proof}

\begin{remark}\label{rem:fact}
Recall from \cite[Sec.~1]{Mum66} that every finite commutative 
group $H$ of order prime to $p$, equipped with a non-degenerate 
alternating pairing $e$, admits a decomposition
$$
H = H_{n_1} \times \cdots \times H_{n_r},
\quad 
e = (e_{d_1}, \ldots, e_{d_r})
$$
such that 
$$
H_{n_i} = \bZ/n_i\bZ \times \cX(\bZ/n_i \bZ) 
\cong (\bZ/n_i\bZ)^2, \quad
e_{d_i}((x,\chi),(x',\chi')) = \chi'(x)^{d_i} \, \chi(x')^{-d_i},
$$
where the $n_i$, $d_j$ are integers satisfying $n_{i+1} \vert n_i$, 
$0 \leq d_i < n_i$, and $(d_i,n_i) = 1$ for all $i$. Moreover,
$n_1, \ldots, n_r$ are uniquely determined by $H$.
Since $H$ is a subgroup of $\wX_n \cong (\bZ/n\bZ)^{2g}$, 
where $g := \dim(X)$, we see that $r \leq g$; conversely, 
any product of $r$ cyclic groups of order prime to $p$ 
can be embedded into $\wX_n$ provided that $r \leq g$.  

It follows that every homogeneous irreducible bundle
admits a decomposition into a product
$$
P = P_1 \cdots P_r,
$$
where each $P_i$ corresponds to $(H_{n_i},e_{d_i})$. Moreover,
the $P_i$ are exactly the irreducible homogeneous bundles
associated with a product of two copies of a cyclic group; we may
call these bundles \emph{cyclic}.

Equivalently, the associated Azumaya algebra satisfies
$$
\cA = \cA_1 \otimes \cdots \otimes \cA_r,
$$
where $\cA_i$ corresponds to $(H_{n_i}, e_{d_i})$.
Moreover, the $\cO_X$-algebra $\cA_i$ is generated by two invertible
sheaves $\cL$, $\cM$ (associated with the natural generators of 
$(\bZ/n_i\bZ)^2$), with relations $x^{n_i} = \xi$, $y^{n_i} = \eta$, 
$xy = \zeta^{d_i} yx$ for any local generators $x \in \cL$, $y \in \cM$,
where $\xi$ (resp.~$\eta$) denotes a local trivialization of 
$\cL^{\otimes n}$ (resp. $\cM^{\otimes n}$), 
and $\zeta$ is a fixed primitive $d_i$th root of unity
(this follows by combining the isomorphism of algebras (\ref{eqn:azu}) 
with the description of the $H_{n_i}$-algebra $\M_{n_i}$ obtained in 
Lemma \ref{lem:cent} (iii)). In particular, $\cA_i$ yields a cyclic 
division algebra over $k(X)$.
\end{remark}

\begin{example}\label{ex:ell}
Let $X$ be an elliptic curve. Then $X$ is canonically isomorphic 
to $\wX$ and the finite subgroups of $X$ admitting a non-degenerate
alternating pairing are exactly the $n$-torsion subgroups $X_n$. 
In view of the above remark, it follows that the irreducible 
homogeneous bundles over $X$ are exactly the cyclic bundles.
By a result of Atiyah (see \cite[Thm.~10]{At57}), they
are exactly the projectivizations of the indecomposable vector 
bundles of coprime rank and degree, i.e., of the simple vector bundles. 
\end{example}

\begin{example}\label{ex:mum}
Returning to an arbitrary abelian variety $X$, we recall from 
\cite[Sec.~1]{Mum66} a geometric construction of Heisenberg groups. 
Let $L$ be a line bundle on $X$, and $K(L)$ the kernel of the 
polarization homomorphism
\begin{equation}\label{eqn:pol}
\varphi_L : X \longrightarrow \wX, \quad x \longmapsto T_x^*(L) \otimes L^{-1}.
\end{equation}
Denoting by $\cG(L)$ the group scheme of automorphisms of the variety $L$ 
which commute with the action of $\bG_m$ by multiplication on fibers,
we have a central extension
$$
1 \longrightarrow \bG_m \longrightarrow \cG(L) \longrightarrow K(L)
\longrightarrow 1.
$$
The associated commutator pairing on $K(L)$ is denoted by $e^L$.

Also, recall that an effective line bundle $L$ is ample if and only if 
$\varphi_L$ is an isogeny; equivalently, $K(L)$ is finite. 
Then the theta group $\cG(L)$ is non-degenerate, and acts on 
the space of global sections $H^0(X,L)$ via its standard representation. 
Thus, $K(L)$ acts on the associated projective space
$$
\vert L \vert := \bP(H^0(X,L))
$$
and the natural map
$$
f: X \times^{K(L)} \vert L \vert \longrightarrow X/K(L) \cong \wX
$$
is an irreducible homogeneous bundle.

As will be shown in detail in Section \ref{sec:exa}, this bundle is 
the projectivization of a natural vector bundle $E$ over $\wX$. 
Moreover, if $X$ is an elliptic curve (so that $X \cong \wX$) 
and $L$ has degree $n$, then $E$ has rank $n$ and degree $-1$.
\end{example}

We now obtain several criteria for a homogeneous projective bundle
to be irreducible:

\begin{proposition}\label{prop:irr}
The following conditions are equivalent for a homogeneous 
bundle $P$:

\smallskip

\noindent
{\rm (i)} $P$ is irreducible.

\smallskip

\noindent
{\rm (ii)} $P$ admits no proper homogeneous sub-bundle.

\smallskip

\noindent
{\rm (iii)} $\ad(P)$ splits into a direct sum of non-zero
algebraically trivial line bundles.

\smallskip

\noindent
{\rm (iv)} $H^0(X,\ad(P)) = 0$.

\smallskip

\noindent
{\rm (v)} $\Aut_X(P)$ is finite.

\medskip

If $P$ is the projectivization of a (semi-homogeneous) vector bundle
$E$, then $P$ is irreducible if and only if $E$ is simple.
\end{proposition}

\begin{proof}
(i)$\Leftrightarrow$(ii) follows from Theorem \ref{thm:str} (iii),
and (i)$\Rightarrow$(iii) from Proposition \ref{prop:irraut} (ii).

(iii)$\Rightarrow$(iv) holds since $H^0(X,L) = 0$ for any
non-zero $L \in \wX$.

(iv)$\Rightarrow$(v) follows from the fact that
$\Lie \, \Aut_X(P) =H^0(X,\ad(P))$.

(v)$\Rightarrow$(i) By Proposition \ref{prop:indec}, $P$ is 
indecomposable and the quotient $\GL(U)^{H_u}/\bG_m \id$ is finite, 
where $U$ is the indecomposable representation of $H_u$ given by that
proposition. But $\GL(U)^{H_u}/\bG_m \id$ has positive dimension for 
any unipotent subgroup scheme $H_u \subset \GL(U)$, unless $\dim(U) = 1$; 
in the latter case, $P$ is clearly irreducible.

The final assertion follows from the equivalence (i)$\Leftrightarrow$(iv)
in view of the isomorphism
$$
H^0(X, \ad(\bP(E))) \cong H^0(X, \End(E))/k \, \id.
$$
\end{proof}

\begin{remark}\label{rem:indec}
The indecomposable homogeneous bundles are exactly 
the products $\bP(U) \, I$, where $U$ is an indecomposable unipotent 
vector bundle, and $I$ an irreducible homogeneous bundle 
(as follows from Proposition \ref{prop:indec}).

In particular, the indecomposable homogeneous bundles over an elliptic 
curve $X$ are exactly the projectivizations $\bP(U \otimes E)$, 
where $U$ is as above, and $E$ is a simple vector bundle (as in 
Example \ref{ex:ell}). 

By a result of Atiyah (see \cite{At57}), any indecomposable vector bundle 
over $X$ is isomorphic to $U \otimes E \otimes L$ for $U$, $E$ as above 
and $L$ a line bundle. Also, $U$ is uniquely determined by its rank;
moreover, $E$ is uniquely determined by its (coprime) rank and degree, 
up to tensoring with a line bundle of degree $0$. 
\end{remark}

Next, we obtain a cohomological criterion for a bundle
to be homogeneous and irreducible, thereby extending a result of 
Mukai about simple semi-homogeneous vector bundles 
(see \cite[Thm.~5.8]{Muk78}): 

\begin{proposition}\label{prop:homirr}
A bundle $P$ is homogeneous and irreducible if and only if 
we have $H^0(X,\ad(P)) = H^1(X,\ad(P)) = 0$; then  
$H^i(X,\ad(P)) = 0$ for all $i \geq 0$.
\end{proposition}

\begin{proof}
Recall that $H^i(X,L) = 0$ for all $i \geq 0$ and all non-zero 
$L \in \wX$. By Proposition \ref{prop:irraut} (ii), the same holds 
with $L$ replaced with $\ad(P)$, if $P$ is homogeneous and irreducible.

For the converse, observe that $\ad(P) = \pi_*(T_{Y/X})^{\PGL_n}$, 
where $\pi: Y \to X$ denotes the $\PGL_n$-torsor associated to $P$,
and $T_{Y/X}$ the relative tangent bundle.
Thus, $\ad(P)$ sits in an exact sequence 
$$
0 \longrightarrow \ad(P) \longrightarrow \pi_*(T_Y)^{\PGL_n} 
\longrightarrow T_X \longrightarrow 0
$$
obtained from the standard exact sequence
$0 \to T_{Y/X} \to T_Y \to \pi^*(T_X) \to 0$ by taking the invariant
direct image under $\pi$. If $H^1(X, \ad(P)) =0$, then the natural map
$$
H^0(Y,T_Y)^{\PGL_n} = H^0(X,\pi_*(T_Y)^{\PGL_n}) \longrightarrow H^0(X,T_X) 
$$
is surjective. But $H^0(Y,T_Y)^{\PGL_n} \cong \Lie(\Aut^{\PGL_n}(Y))$ 
and $H^0(X,T_X) \cong \Lie(\Aut(X))$; moreover, 
$\Aut^{\PGL_n}(Y) = \Aut(P)$ is smooth by Remark \ref{rem:form},
and $\Aut(X)$ is smooth as well. Hence the homomorphism 
$\Aut^{\PGL_n}(Y) \to \Aut(X)$
is surjective on neutral components, i.e., $Y$ is homogeneous. Thus,
$P$ is homogeneous, too. If in addition $H^0(X,\ad(P)) = 0$, 
then $P$ is irreducible by Proposition \ref{prop:irr}.
\end{proof}

\begin{remark}\label{rem:coh}
The above argument shows that a bundle $P$ is homogeneous if it satisfies
$H^1(X,\ad(P)) = 0$. This may also be seen as follows:
observe that $\ad(P) = f_*(T_{P/X})$ (as follows e.g. by considering
an \'etale trivialization of $P$). Moreover, $R^if_*(T_{P/X}) = 0$ 
for all $i \geq 1$, since $H^i(\bP^{n-1},T_{\bP^{n-1}}) = 0$ 
for all such $i$. As a consequence, $H^1(P,T_{P/X}) = 0$. 
Then $f$ is rigid as a morphism with target $X$ 
in view of \cite[Cor.~3.4.9]{S06}. It follows readily that $P$ is
homogeneous.

The converse statement does not hold, e.g., when $X$ is an elliptic curve
in characteristic zero, $U_n$ is the indecomposable unipotent vector bundle 
of rank $n \geq 2$, and $P = \bP(U_n)$. Then 
$$
\ad(P) \cong (U_n \otimes U_n^*)/k \, \id \cong
U_{2n-1} \oplus U_{2n-3} \oplus \cdots \oplus U_3
$$
and hence $H^0(X,\ad(P))$ has dimension $n -1$. By the 
Riemann-Roch theorem, the same holds for $H^1(X,\ad(P))$.
\end{remark}

\subsection{Projectivizations of vector bundles}
\label{subsec:pvb}

In this subsection, we characterize those homogeneous projective
bundles that are projectivizations of (not necessarily homogeneous)
vector bundles. We first consider a special class of bundles, defined
as follows. 

Given a positive integer $m$, not divisible by $p$,
we say that a bundle $P$ is \emph{trivialized by $m_X$} 
(the multiplication by $m$ in $X$) if the pull-back bundle 
$m_X^*(P) \to X$ is trivial. 

In fact, every such bundle is homogeneous, as a consequence of
the following:

\begin{proposition}\label{prop:stir}
{\rm (i)} A bundle $P$ is trivialized by $m_X$ 
if and only if $P \cong X \times^{X_m} \bP^{n-1}$ as bundles over 
$X \cong X/X_m$, for some action of $X_m$ on $\bP^{n-1}$. 

\smallskip

\noindent
{\rm (ii)} Any irreducible homogeneous $\bP^{n-1}$-bundle
is trivialized by $n_X$.
\end{proposition}

\begin{proof}
(i) If $P$ is trivialized by $m_X$, then we have a cartesian square
$$
\CD
X \times \bP^{n-1} @>{p_1}>> X \\
@V{q}VV @V{m_X}VV \\
P @>{f}>> X, \\
\endCD
$$
where $p_1$ denotes the first projection.
Thus, the action of $X_m$ by translations on $X$ lifts to an action 
on $X \times \bP^{n - 1}$ such that $q$ is invariant. This action is of 
the form 
$$
x \cdot (y,z) = (x + y, \varphi(x,y) \cdot z)
$$
for some morphism $\varphi : X_m \times X \to \Aut(\bP^{n-1}) = \PGL_n$.
But every morphism from the abelian variety $X$ to the affine variety 
$\PGL_n$ is constant. Thus, $\varphi$ is independent of $y$, i.e., 
$\varphi$ yields an action of $X_m$ on $\bP^{n-1}$. Moreover, the 
$X_m$-invariant morphism $q$ factors through a morphism of 
$\bP^{n-1}$-bundles $X \times^{X_m} \bP^{n-1} \to P$ 
which is the desired isomorphism.

The converse implication is obvious.

(ii) Write $P = G \times^H \bP^{n-1}$ as in Theorem \ref{thm:str};
then $H$ is killed by $n$ in view of the structure of non-degenerate 
theta groups. In other words, the homomorphism $\gamma : G \to X$ 
is an isogeny with kernel killed by $n$. Thus, there exists a unique 
isogeny $\tau : X \to G$ such that $\gamma \tau = n_X$. 
Then $X_n = \tau^{-1}(H)$ and hence $X = X \times^{X_n} \bP^{n-1}$, 
where $X_n$ acts on $\bP^{n-1}$ via the surjective homomorphism 
$\tau_{\vert X_n} : X_n \to H$.
\end{proof}

By the above proposition, a bundle $P$ trivialized by $m_X$ 
defines an alternating bilinear map
$$
e_{P,m} : X_m \times X_m \longrightarrow \mu_m.
$$
Moreover, the irreducible homogeneous bundles are 
classified by those maps such that $[X_m : X_m^{\perp}] = m^2$ 
(as follows from Proposition \ref{prop:irr}).
Also, one easily checks that the assignement 
$P \mapsto e_{P,m}$ is multiplicative, i.e.,
$e_{P_1 P_2,m} = e_{P_1,m} \, e_{P_2,m}$ and $e_{P^*,m} = e_{P,m}^{-1}$.

We may now obtain the desired characterization:

\begin{proposition}\label{prop:semi}
Let $P$ be a bundle trivialized by $m_X$. Then $P$ is the projectivization 
of a vector bundle if and only if there exists a line bundle $L$ on $X$ 
such that $e_{P,m} = e^{L^{\otimes m}}\vert_{X_m}$ (this makes sense as
$K(L^{\otimes m})$ contains $X_m$). 
\end{proposition}

\begin{proof}
Assume that $P = \bP(E)$ for some vector bundle $E$ 
of rank $n$ on $X$. Since the projective bundle $m_X^* (\bP(E))$ 
is trivial, we have 
$$
m_X^*(E) \cong M^{\oplus n}
$$ 
for some line bundle $M$ on $X$. 
Replacing $E$ with $E \otimes N$, where $N$ is a symmetric 
line bundle on $X$, leaves $\bP(E)$ unchanged and replaces 
$m_X^*(E)$ with $m_X^*(E) \otimes N^{\otimes m^2}$, and hence $M$ 
with $M \otimes N^{\otimes m^2}$. Taking for $N$ a large power of an
ample symmetric line bundle, we may assume that $M$ is very ample. 

The pull-back $m_X^*(E)$ is equipped with an $X_m$-linearization.
Equivalently, the action of $X_m$ by translations on $X$ lifts to
an action on $M^{\oplus n}$ which is linear on fibers. In particular,
$T_x^*(M^{\oplus n}) \cong M^{\oplus n}$ for any $x \in X_m$. 
This isomorphism is given by an $n \times n$ matrix of maps 
$T_x^*M \to M$; thus, $H^0(X,T_x^*(M^{-1}) \otimes M) \neq 0$.
Since $T_x^*(M^{-1}) \otimes M \in \wX$, it follows that this line
bundle is trivial. In other words, $X_m \subset K(M)$; this is 
equivalent to the existence of a line bundle $L$ in $X$ such that
$M = L^{\otimes m}$. Moreover, we have a representation of $X_m$ in 
$H^0(X, M^{\oplus n}) \cong H^0(X,M) \otimes k^n$
that lifts the homomorphism 
\begin{equation}\label{eqn:phi}
\phi : X_m \longrightarrow \PGL(H^0(X,M)) \times \PGL_n
\end{equation}
given by the $X_m$-action on $\bP(H^0(X,M))$ as a subgroup of $K(M)$,
and the $X_m$-action on $\bP^{n-1}$ that defines $P$. It follows that 
$e^M \, e_{P,m} = 1$ on $X_m$; equivalently, $e_{P,m}$ is the restriction 
to $X_m$ of 
$e^{M^{\otimes(-1)}} = e^{L^{\otimes(-m)}} = e^{L^{\otimes m(m-1)}}$
(since $e^{L^{\otimes m^2}} = 1$).

To show the converse, we reduce by inverting the above arguments
to the case that $e^M \, e_{P,m} = 1$ on $X_m$ for some line bundle 
$M$ on $X$ such that $X_m \subset K(M)$; we may also assume that 
$M$ is very ample. Then $X_m$ acts on $H^0(X,M^{\oplus n})$ by lifting 
the homomorphism (\ref{eqn:phi}). Moreover, the evaluation morphism
$$
\cO_X \otimes H^0(X,M^{\oplus n}) = \cO_X \otimes H^0(X,M) \otimes k^n
\longrightarrow M \otimes k^n = M^{\oplus n}
$$
is surjective and its kernel is stable under the induced action 
of $X_m$ (since the analogous morphism $\cO_X \otimes H^0(X,M) \to M$
is equivariant with respect to the theta group of $X_m \subset K(M)$). 
Thus, $X_m$ acts on $M^{\oplus n}$ by lifting its action on $X$ via 
translation. Now $M^{\oplus n}$ descends to a vector bundle on 
$X/X_m \cong X$ with projectivization $P$.
\end{proof}

Next, we extend the statement of Proposition \ref{prop:semi} to all 
homogeneous bundles $P$. We use the notation of Section \ref{sec:str}; 
in particular, the associated pairing $e_F$ introduced in Proposition 
\ref{prop:theta}. Then $e_F$ factors through a non-degenerate pairing 
on $F/F^{\perp}\cong H/H^{\perp}$ and this group is killed by the 
homogeneous index $d = d(H)$ defined by (\ref{eqn:hindex}). 
Thus, the isogeny $G/H^{\perp} \to G/H = X$ has its kernel killed by 
$d$; as in the proof of Proposition \ref{prop:stir} (ii), 
this yields a canonical surjective homomorphism $X_d \to H/H^{\perp}$ 
and, in turn, a bilinear alternating pairing $e_P$ on~$X_d$.

\begin{theorem}\label{thm:obs}
With the above notation, $P$ is the projectivization of a vector bundle 
if and only if $e_P = e^{L^d}\vert_{X_d}$ for some line bundle $L$ on $X$.
\end{theorem}

\begin{proof}
Choose a linear subspace $S\subset \bP^{n-1}$ which is $H$-stable, 
and minimal for this property. Then $S$ yields a homogeneous irreducible 
$\bP^{d-1}$-sub-bundle of $P$ and the associated pairing on $X_d$ 
is just $e_P$. Now the statement is a consequence of Proposition 
\ref{prop:semi} together with the following observation.

\begin{lemma}\label{lem:sub}
Let $f :P \to Z$ be a projective bundle over a non-singular variety, and
$f_1: P_1 \to Z$ a projective sub-bundle. Then $P$ is the projectivization
of a vector bundle if and only if so is $P_1$.
\end{lemma}

\begin{proof}
Clearly, if $P = \bP(E)$ for some vector bundle $E$ over $Z$, then 
$P_1 = \bP(E_1)$ for some sub-bundle $E_1 \subset E$. To show the converse,
consider the $\PGL_n$-torsor $\pi: Y \to Z$ associated with $P$; 
recall that the sub-bundle $P_1$ yields a reduction of structure group 
to a $\PGL_{n,n_1}$-torsor $\pi_1 : Y_1 \to Z$, where 
$\PGL_{n,n_1}\subset \PGL_n$ denotes the stabilizer of 
$\bP^{n_1 -1} \subset \bP^{n - 1}$. 
We have an exact sequence of algebraic groups
$$
\CD
1 @>>> G_{n,n_1} @>>> \PGL_{n,n_1} @>{r}>> \PGL_{n_1} @>>> 1,
\endCD
$$
where $r$ denotes the restriction to $\bP^{n_1 -1}$ and 
$G_{n,n_1} \cong \M_{n_1,n-n_1} \ltimes \GL_{n-n_1}$, 
the semi-direct product being defined by the natural action of 
$\GL_{n-n_1}$ on the space of matrices $M_{n_1,n-n_1}$. 
Also, $\pi_1$ factors as 
$$
\CD
Y_1 @>{\varphi}>> Y_1/G_{n,n_1} @>{\psi}>> Z,
\endCD
$$
where $\varphi$ is a $G_{n,n_1}$-torsor and $\psi$ is the $\PGL_{n_1}$-torsor 
associated with $P_1$. By assumption, $P_1 = \bP(E_1)$ for some vector bundle 
$E_1$; this is equivalent to $\psi$ being locally trivial in view of 
\cite[Prop.~18]{Se58}. But $\varphi$ is locally trivial as well, 
since the algebraic group $G_{n,n_1}$ is special by 
[loc.~cit., 4.3, 4.4]. Thus, $\pi_1$ is locally trivial, and hence so is 
$\pi$. We conclude that $P = \bP(E)$ for some vector bundle $E$.

Alternatively, one may use the fact that $P$ is the projectivization
of a vector bundle if and only if $f$ has a rational section 
(\cite[Prop.~18]{Se58} again), and conclude by applying
\cite[Prop.~5.3.1]{GS06}. 
\end{proof}

\end{proof}

\begin{remark}\label{rem:ber}
We now relate Proposition \ref{prop:semi} to a description of 
the Brauer group $\Br(X)$, due to Berkovich. 
Recall from \cite[Sec.~I.8.4]{Gr68} that $\Br(X)$ may be viewed as
the set of equivalence classes of projective bundles over $X$, where
two such bundles $P_1$, $P_2$ are equivalent if there exist vector bundles
$E_1$, $E_2$ such that $\bP(E_1) P_1 \cong \bP(E_2) P_2$; the group 
structure stems from the operations of product and duality. By 
\cite[Sec.~3]{Be72}, we have an exact sequence for any positive 
integer $n$:
$$
\CD
0 @>>> \Pic(X)/n \Pic(X) @>{\varphi}>> \Hom(\Lambda^2 X_n,\mu_n) 
@>{\psi}>> \Br(X)_n @>>> 0,
\endCD
$$
where $\Hom(\Lambda^2 X_n,\mu_n)$ consists of the bilinear alternating
pairings $X_n \times X_n \to \mu_n$ and $\Br(X)_n \subset \Br(X)$ 
denotes the $n$-torsion subgroup; the map $\varphi$ sends the class of 
$L \in \Pic(X)$ to the pairing $e^{L^{\otimes n}}\vert_{X_n}$ and the map
$\psi$ sends $e$ to the class of the Azumaya algebra
$$
\cA := \bigoplus_{\alpha \in \wX_n, \, \sigma \in X_n}
\cL_{\alpha} \, e_{\sigma},
$$
where $\cL_{\alpha}$ denotes the invertible sheaf associated with
$\alpha$ and the multiplication is defined by
$$
f_{\alpha} \, e_{\sigma} \cdot f_{\beta} \, e_{\tau} 
= \bar{e}_n(\beta,\sigma) \, a_{\sigma,\tau} \, f_{\alpha} \, f_{\beta}
\, e_{\sigma + \tau}.
$$
Here $f_{\alpha}$ (resp. $f_{\beta}$) is a local section of 
$\cL_{\alpha}$ (resp. $\cL_{\beta}$); $\bar{e}_n$ is the canonical 
pairing between $\wX_n$ and $X_n$, and 
$\{ a_{\sigma,\tau} \} \in Z^2(X_n,\bG_m)$
is a $2$-cocycle such that 
$e(\sigma,\tau) = a_{\sigma,\tau} \, a_{\tau,\sigma}^{-1}$.
(The class of $\cA$ in the Brauer group does not depend on the choice
of the representative $\{a_{\sigma,\tau} \}$ of $e$ viewed as an element
of $H^2(X_n,\bG_m)$.) Thus,
$$
\cL := \bigoplus_{\alpha \in \wX_n} \cL_{\alpha} \, e_0
$$
is a maximal \'etale subalgebra of $\cA$ in the sense of 
\cite[D\'ef.~5.6]{Gr68}; note that $\cL \cong (n_X)_* \cO_X$
as $\cO_X$-algebras. Moreover, the left $\cL$-module $\cA$ is free
with basis $(a_{\sigma})_{\sigma \in X_n}$. By 
[loc.~cit., Cor.~5.5], it follows that 
$n_X^*(\cA) \cong \M_m(\cO_X)$, where $m := \#(X_n) = n^{2g}$.
In other words, the projective bundle associated with $\cA$ is 
trivialized by $n_X$. In view of Proposition \ref{prop:stir}, it follows
that the associated projective bundle is homogeneous.

In fact, \emph{any class in $\Br(X)_n$ is represented by an irreducible
homogeneous bundle}. Indeed, given any homogeneous bundle $P$, we may
choose an irreducible sub-bundle $P_1$; then the product $P_1 P_1^*$ is a
sub-bundle of $P P_1^*$ and is the projectivization of a vector bundle.
By Lemma \ref{lem:sub}, it follows that the class of $P P_1^*$ in $\Br(X)$
is trivial; equivalently, $P$ and $P_1$ have the same class there.

Also, recall that the natural map $\Br(X) \to \Br(k(X))$ is injective
(see \cite[Sec.~II.1]{Gr68}).
As a very special case of a theorem of Merkurjev and Suslin 
(see \cite[Thm.~2.5.7]{GS06}), each class in $\Br(k(X))_n$ can be 
represented by a tensor product of cyclic algebras. So the decomposition
of classes in $\Br(X)_n$ obtained in Remark \ref{rem:fact} may be viewed
as a global analogue of that result for abelian varieties.

Finally, note that Proposition \ref{prop:semi} is equivalent to 
the assertion that 
\emph{the image of $\varphi$ consists of those pairings associated 
with projectivizations of semi-homogeneous vector bundles.} 
In loose terms, the Brauer group is generated by homogeneous bundles
and the relations arise from semi-homogeneous vector bundles.
\end{remark}

\section{Examples}
\label{sec:exa}

Let $X$ be an abelian variety, and $\lambda$ an effective class 
in the N\'eron-Severi group $NS(X)$ viewed as the group of divisors 
on $X$ modulo algebraic equivalence. 
The effective divisors on $X$ with class $\lambda$ are parametrized 
by a projective scheme $\Div^{\lambda}(X)$. Indeed, the Hilbert polynomial 
of any such divisor $D$, relative to a fixed ample line bundle on $X$, 
depends only on $\lambda$; thus, $\Div^{\lambda}(X)$ is a union 
of connected components of the Hilbert scheme $\Hilb(X)$. 

Also, recall that the line bundles on $X$ with class $\lambda$ are 
parametrized by the Picard variety $\Pic^{\lambda}(X)$. Choosing
$L$ in that variety, we have 
$$
\Pic^{\lambda}(X) = L \otimes \Pic^0(X) = L \otimes \wX.
$$
On $X \times \Pic^{\lambda}(X)$ we have a universal bundle: 
the Poincar\'e bundle $\cP$, uniquely determined up to the pull-back 
of a line bundle under the second projection
$$
\pi : X \times \Pic^{\lambda}(X) \longrightarrow \Pic^{\lambda}(X).
$$
The universal family on $\Div^{\lambda}(X)$ yields a morphism
\begin{equation}\label{eqn:div}
f : \Div^{\lambda}(X) \longrightarrow \Pic^{\lambda}(X), 
\quad D \longmapsto \cO_X(D).
\end{equation}
Note that $X$ acts on $\Div^{\lambda}(X)$ and on $\Pic^{\lambda}(X)$ 
via its action on itself by translations; moreover, $f$ is equivariant. 
Also, the isotropy subgroup scheme in $X$ of any point of $\Pic^{\lambda}(X)$
is the group scheme $K(L)$ that occured in Example \ref{ex:mum}. 

If $\lambda$ is ample, then $\Pic^{\lambda}(X)$ is the $X$-orbit
$X \cdot L \cong X/K(L)$. 
Thus, $f$ is a homogeneous fiber bundle over $X/K(L)$; the latter 
abelian variety is isomorphic to $\wX$ via the polarization homomorphism
(\ref{eqn:pol}). 

\begin{proposition}\label{prop:exa}
Let $\lambda \in \NS(X)$ be an ample class, and 
$L \in \Pic^{\lambda}(X)$.

\smallskip

\noindent 
{\rm (i)} We have an isomorphism
$$
\Div^{\lambda}(X) \cong X \times^{K(L)} \vert L \vert
$$
of homogeneous bundles over $X/K(L)$. In particular, $\Div^{\lambda}(X)$ 
is a homogeneous projective bundle over $\wX$.

\smallskip

\noindent
{\rm (ii)} The sheaf $\cE := \pi_*(\cP)$ is locally free, and 
the morphism (\ref{eqn:div}) is the projectivization of the 
corresponding vector bundle. 

\smallskip

\noindent
{\rm (iii)} The group scheme $\Aut(\Div^{\lambda}(X))$ 
is the semi-direct product of $X$ (acting by translations) 
with the subgroup of $\Aut_{\gp}(X)$ that preserves $K(L)$ and $e^L$.
\end{proposition} 

\begin{proof}
(i) Clearly, the set-theoretic fiber of $f$ at $L$ is the projective
space $\vert L \vert$, and its dimension $h^0(L) -1 = \chi(L) -1$ is
independent of $L \in \Pic^{\lambda}(X)$. As a consequence, the scheme
$\Div^{\lambda}(X)$ is irreducible of dimension $\dim(X) + h^0(L) -1$. 

To complete the proof, it suffices to show that the differential of $f$ 
at any $D \in \vert L \vert$ is surjective with kernel of dimension 
$h^0(L) -1$. Identifying $\Div^{\lambda}(X)$ with a union of components 
of $\Hilb(X)$, and $\Pic^{\lambda}(X)$ with $\wX$, the differential
$$
T_D f : T_D \Div^{\lambda}(X) \longrightarrow T_L \Pic^{\lambda}(X)
$$
is identified with the boundary map 
$\partial : H^0(D, L_{\vert D}) \to H^1(X,\cO_X)$
of the long exact sequence of cohomology associated with 
the short exact sequence
$$
0 \longrightarrow \cO_X \longrightarrow L \longrightarrow 
L_{\vert D} \longrightarrow 0
$$ 
(see \cite[Prop.~3.3.6]{S06}).
Since  $H^1(X,L) = 0$, this long exact sequence begins with
$$
\CD
0 @>>> k @>>> H^0(X,L) @>>> H^0(D,L_{\vert D}) 
@>{\partial}>> H^1(X,\cO_X) @>>> 0
\endCD 
$$
which yields the desired assertion.

(ii) The vanishing of $H^1(X,L)$ also implies that $\cE$ is locally free
and satisfies $\cE(L) \cong H^0(X,L)$. Thus, it suffices to check that 
the associated projective bundle $\bP(\cE)$ is homogeneous.
But for any $x \in X$, there exists an invertible sheaf $L_x$ on 
$\Pic^{\lambda}(X)$ such that
$$
(T_x,T_x)^* (\cP) \cong \cP \otimes \pi^* L_x
$$ 
in view of the universal property of the Poincar\'e bundle $\cP$.
Since $\pi_*(T_x,T_x)^*(\cP) \cong T_x^* (\pi_*(\cP)) = T_x^*(\cE)$, 
this yields an isomorphism
$$
T_x^* (\cE) \cong \cE \otimes L_x.
$$
In other words, $\cE$ is semi-homogeneous.

(iii) is checked by arguing as in the proof of Proposition 
\ref{prop:irraut} (iii).
\end{proof}

The case of an arbitrary effective class $\lambda$ reduces to 
the ample case in view of the following:

\begin{proposition}\label{exe}
Let $\lambda \in \NS(X)$ be an effective class, $L \in \Pic^{\lambda}(X)$, 
and $q : X \to \bar{X}$ the quotient map by the reduced neutral component 
$K(L)^0_{\red} \subset K(L)$. Then $\lambda = q^*(\bar{\lambda})$ for 
a unique ample class $\bar{\lambda} \in \NS(\bar{X})$, and 
$f : \Div^{\lambda}(X) \to \Pic^{\lambda}(X)$ may be identified with 
$\bar{f} : \Div^{\bar{\lambda}}(\bar{X}) \to \Pic^{\bar{\lambda}}(\bar{X})$.
\end{proposition}

\begin{proof}
We claim that any $D \in \Div^{\lambda}(X)$ equals $q^*(\bar{D})$ 
for some ample effective divisor $\bar{D}$ on $\bar{X}$.

To see this, recall that $n D$ is base-point-free for any $n \geq 2$; 
this yields morphisms
$$
\gamma_n : X \longrightarrow \bP(H^0(X, L^{\otimes n})^*) \quad (n \geq 2)
$$
which are equivariant for the action of $K(L)$. The abelian variety
$K(L)^0_{\red}$  acts trivially on each projective space 
$\bP(H^0(X, L^{\otimes n})^*)$; thus, each $\gamma_n$ is invariant 
under $K(L)^0_{\red}$. In the Stein factorization of $\gamma_n$ as
$$
\CD 
X @>{\varphi_n}>> Y_n @>{\psi_n}>> \bP(H^0(X,L^{\otimes n})^*),
\endCD
$$
where $(\varphi_n)_*(\cO_X) = \cO_{Y_n}$ and $\psi_n$ is finite, 
the morphism $\varphi_n$ is the natural map
$$
\varphi : 
X \longrightarrow \Proj \bigoplus_{m=0}^{\infty} H^0(X,L^{\otimes m})
=: Y.
$$
In particular, $\varphi_n$ is independent of $n$ and invariant under 
$K(L)^0_{\red}$. Moreover, since $nD$ is the pull-back of a hyperplane under 
$\gamma_n$ for any $n \geq 2$, we see that $D = 3D - 2D = \varphi^*(E)$ 
for some Cartier divisor $E$ on $Y$. Then $E$ is effective and 
$H^0(X,L^{\otimes n}) \cong H^0(Y,M^{\otimes n})$ for all $n$,
where $M := \cO_Y(E)$; it follows that $E$ is ample. 
Consider the factorization
$$
\bar{\varphi} : \bar{X} := X/K(L)^0_{\red} \longrightarrow Y,
$$
the effective divisor $\bar{D}:= \bar{\varphi}^*(E)$, and the associated
invertible sheaf $\bar{L} = \bar{\varphi}^*(M)$. Then 
$L = q^*(\bar{L})$. Thus, the group scheme $K(\bar{L}) = K(L)/K(L)^0_{\red}$ 
is finite and $\bar{L}$ has non-zero global sections; hence $\bar{L}$ is 
ample. Thus, $\bar{\varphi}$ is finite. But
$\bar{\varphi}_*(\cO_{\bar{X}}) = \cO_Y$; it follows that $\bar{\varphi}$ 
is an isomorphism, and this identifies $\varphi$ with $q$. 
This proves the claim.

As a consequence, $\lambda = q^*(\bar{\lambda})$ for a unique 
ample class $\bar{\lambda}$. We now show that the morphism 
$$ 
q^* : \Div^{\bar{\lambda}}(\bar{X}) \longrightarrow \Div^{\lambda}(X)
$$
is an isomorphism. By the first step, $q^*$ is bijective. 
In view of Proposition \ref{prop:exa}, it follows that the scheme
$\Div^{\lambda}(X)$ is irreducible of dimension 
$\dim(\bar{X}) + h^0(\bar{X},\bar{L}) - 1$. On the other hand, 
the Zariski tangent space of $\Div^{\lambda}(X)$ at $D$ equals
$$
H^0(D,L_{\vert D}) \cong H^0(\bar{D}, \bar{L}_{\vert \bar{D}}) 
= T_{\bar{D}} \Div^{\bar{\lambda}}(\bar{X}).
$$
Thus, $q^*$ is \'etale and hence is an isomorphism.
\end{proof}

In the above construction, one may replace the abelian variety $X$ with 
any smooth projective variety; for example, a curve $C$. Then an effective 
class in $\NS(C) \cong \bZ$ is just a non-negative integer $d$. Moreover, 
$\Div^d(C)$ is the symmetric product $C^{(d)}$, a smooth projective
variety of dimension $d$ equipped with a morphism
\begin{equation}\label{eqn:cur}
f = f_d: C^{(d)} \longrightarrow \Pic^d(C).
\end{equation}
Choosing a point of $C$, we may identify $\Pic^d(C)$ with the
Jacobian variety $J = J(C)$.

If $d > 2g - 2$, where $g$ denotes of course the genus of $C$, 
then $f$ is the projectivization of a vector bundle $E = E_d$ 
on $\Pic^d(C)$, the direct image of the Poincar\'e bundle on 
$C \times \Pic^d(C)$ under the second projection. Moreover,
$E$ has rank $n := d - g + 1$.

\begin{proposition}\label{prop:nex}
With the above notation, the projective bundle (\ref{eqn:cur}) 
is homogeneous if and only if $g \leq 1$.
\end{proposition}

\begin{proof}
Assume that (\ref{eqn:cur}) is homogeneous. Then $E$ is semi-homogeneous; 
in view of \cite[Lem.~6.11]{Muk78}, we then have an isomorphism
of vector bundles on $J$
$$
n_J^*(E) \cong \det(E)^{\otimes n} \otimes F
$$
for some homogeneous vector bundle $F$.    
Moreover, the Chern classes of $F$ are algebraically trivial by 
\cite[Thm.~4.17]{Muk78}. Thus, the total Chern class of $E$
satisfies 
$$
n_J^*(c(E)) = (1  + n c_1(E))^n
$$ 
in the cycle ring of $J$ modulo algebraic equivalence. Since 
$n_J^*(c_1(E)) = n^2 c_1(E)$ 
in that ring, this yields
\begin{equation}\label{eqn:chern}
c(E) = \big( 1 + \frac{c_1(E)}{n} \big)^n .
\end{equation}

We now recall a formula for $c(E)$ due to Mattuck (see 
\cite[Thm.~3]{Ma61}). Denoting by $W_i$ the image of $p_i$ for 
$0 \leq i \leq g$, we have
$$
c(E) = \sum_{i=0}^g  (-1)^i [W_{g-i}^-],
$$ 
where $W_j^-$ denotes the image of $W_j$ under the involution 
$(-1)_J$ and the equality holds again modulo algebraic 
equivalence. In particular, 
$$
c_1(E) = - [W_{g-1}^-] = - \theta,
$$ 
where $\theta$ denotes the Chern class of the theta divisor, and 
$$
c_g(E) = (-1)^g e,
$$ 
where $e$ denotes the class of a point. In view of (\ref{eqn:chern}), 
this yields
$$
e = {n \choose g} \frac{{\theta}^g}{n^g} \, .
$$ 
Since $\theta^g = g! \, e$, we obtain $n^g = n(n-1) \cdots (n - g + 1)$
and hence $g \leq 1$.

Conversely, if $g = 0$ then $C^{(d)} = \bP^d$ and there is nothing
to prove; if $g = 1$ then the assertion follows from Proposition 
\ref{prop:exa}.
\end{proof}

\begin{remark}\label{rem:ex}
By \cite{EL92}, the vector bundle $E$ is stable with respect to 
the principal polarization of $J$. In particular, $E$ is simple,
i.e., $\Aut_J(P)$ is finite. This yields examples of 
simple vector bundles on abelian varieties which are not 
semi-homogeneous (see \cite{Od71} for the first construction of 
bundles satisfying these properties).
\end{remark}

\section{Homogeneous self-dual projective bundles}
\label{sec:hsdpb}

\subsection{Generalities on self-dual bundles}
\label{subsec:gensdpb}

Throughout this subsection, we assume that $p \neq 2$;
we consider projective bundles over a fixed variety $X$. 
Let $f: P \to X$ be a $\bP^{n-1}$-bundle, and $f^* : P^* \to X$ 
the dual bundle. By contravariance, any isomorphism of bundles 
\begin{equation}\label{eqn:iso}
\varphi: P \longrightarrow P^*
\end{equation}
defines a dual isomorphism $\varphi^*: P = P^{**} \to P^*$.
We say that (\ref{eqn:iso}) is \emph{self-dual} if 
$\varphi^* = \varphi$. 

For later use, we now present some general results on self-dual 
bundles; we omit their (easy) proofs, which can be found 
in the arXiv version of this article \cite{Brlong}.

\begin{proposition}\label{prop:sd}
Given a $\bP^{n-1}$-bundle $P$, there is a bijective correspondence 
between the self-dual morphisms (\ref{eqn:iso}) and the reductions of 
structure group of the associated $\PGL_n$-torsor $\pi: Y \to X$ 
to a $\PO_{n,\varepsilon}$-torsor $\psi: Z \to X$, where 
$\varepsilon = \pm 1$ and $\PO_{n,\varepsilon} \subset \PGL_n$ 
denotes the projective orthogonal (resp.~symplectic) group 
if $\varepsilon = +1$ (resp.~$-1$).
\end{proposition}

We say that the self-dual morphism (\ref{eqn:iso}) is 
\emph{symmetric} (resp. \emph{alternating}) if $\varepsilon = 1$
(resp. $= -1$). Denote by $\GO_{n,\varepsilon}$ the preimage of
$\PO_{n,\varepsilon}$ in $\GL_n$. Then $\GO_{n,\varepsilon}$
is the stabilizer of a unique line in the space of bilinear forms
on $k^n$. Moreover, any such semi-invariant form $B$ is 
non-degenerate, and symmetric (resp. alternating) if so is 
$\varphi$.

The group $\GO_{n,\varepsilon}$ is connected and reductive for any 
$n$; hence so is $\PO_{n,\varepsilon}$. If $n$ is odd, then we must
have $\varepsilon = + 1$, and $\PO_{n,\varepsilon} = \SO_n$; if $n$
is even, then $\PO_{n, +1} = \PSO_n$ and $\PO_{n,-1} = \PSp_n$.
As a consequence, $\PO_{n,\varepsilon}$ is semi-simple of adjoint
type unless $n = 2$ and $\varepsilon = 1$; then 
$\PO_{2, + 1} = \bG_m$.

Together with the results of \cite{Gr68} recalled in Subsection
\ref{subsec:genpb}, Proposition \ref{prop:sd} yields one-to-one
correspondences between \emph{self-dual $\bP^{n-1}$-bundles} 
(i.e., bundles equipped with a self-dual morphism), 
$\PO_{n,\varepsilon}$-torsors, and Azumaya algebras $\cA$ of rank 
$n^2$ equipped with an involution (as in \cite{PS92}); these 
correspondences preserve morphisms. The $\PO_{n,\varepsilon}$-torsor
$Z \to X$ corresponds to the  associated bundle 
$P = Z \times^{\PO_{n,\varepsilon}} \bP^{n-1} \to X$ equipped with
the isomorphism to $P^*$ arising from the 
$\PO_{n,\varepsilon}$-equivariant isomorphism 
$\bP^{n-1} \stackrel{\cong}{\longrightarrow} (\bP^{n-1})^*$
given by $B$. The associated Azumaya algebra is 
the sheaf of local sections of the matrix bundle
$Z \times^{\PO_{n,\varepsilon}} \M_n$ equipped with the involution
arising from the isomorphism $\M_n \to (\M_n)^{\op}$ defined by
the adjoint with respect to the pairing $B$. 

Like for $\bP^{n-1}$-bundles, we may define the 
\emph{product} of the self-dual bundles $(P_i,\varphi_i)$ 
($i = 1, 2$) in terms of the associated 
$\PO_{n_i,\varepsilon_i}$-torsors $Z_i \to X$. Specifically,
the product $(P_1P_2,\varphi_1\varphi_2)$ corresponds to the 
$\PO_{n_1 n_2, \varepsilon_1 \varepsilon_2}$-torsor obtained from 
the $\PO_{n_1,\varepsilon_1} \times \PO_{n_2,\varepsilon_2}$-torsor
$Z_1 \times_X Z_2 \to X$ by the extension of structure groups
$$
\PO_{n_1,\varepsilon_1} \times \PO_{n_2,\varepsilon_2}
= \PO_{\varepsilon_1}(k^{n_1}) \times \PO_{\varepsilon_2}(k^{n_2})
\stackrel{\rho}{\longrightarrow}
\PO_{\varepsilon_1 \varepsilon_2}(k^{n_1} \otimes k^{n_2}) 
= \PO_{n_1 n_2, \varepsilon_1 \varepsilon_2},
$$
where $\rho$ stems from the natural map
$\GO_{\varepsilon_1}(k^{n_1}) \times \GO_{\varepsilon_2}(k^{n_2})
\to \GO_{\varepsilon_1 \varepsilon_2}(k^{n_1} \otimes k^{n_2})$.
This product also corresponds to the tensor product of 
algebras with involutions, as considered in \cite{PS92}. 

Next, we introduce a notion of \emph{decomposition} of self-dual bundles;
for this, we need some observations on duality for sub-bundles.
Any $\bP^{n_1 - 1}$-sub-bundle $P_1$ of a bundle $P$ defines 
a $\bP^{n - n_1 - 1}$-sub-bundle of $P^*$, as follows:
$P_1$ corresponds to a $\PGL_n$-equivariant morphism $\gamma$
from $Y$ to the Grassmannian $\PGL_n/\PGL_{n,n_1}$ and hence
to an equivariant morphism $\gamma^*$ from $Y^*$ to the dual 
Grassmannian, $\PGL_n/\PGL_{n, n - n_1}$. The latter morphism 
yields the desired sub-bundle $P_1^{\perp}$. One checks that
$P_1^{\perp \perp} = P_1$ under the identification of $P$ with $P^{**}$.
Moreover, every decomposition $(P_1,P_2)$ of $P$ yields a 
decomposition $(P_2^{\perp},P_1^{\perp})$ of $P^*$, of the same type.  
We may now define a \emph{decomposition} of a self-dual bundle 
$(P,\varphi)$ as a decomposition $(P_1,P_2)$ of the bundle $P$, 
such that $\varphi(P_1) = P_2^{\perp}$; then also 
$\varphi(P_2) = P_1^{\perp}$ by self-duality.

\begin{proposition}\label{prop:dec}
Under the corrrespondence of Proposition \ref{prop:sd}, 
the decompositions of type $(n_1,n_2)$ of $(P,\varphi)$ correspond
bijectively to the reductions of structure group of the 
$\PO_{n,\varepsilon}$-torsor $Z$ to a 
$\P(\O_{n_1,\varepsilon} \times \O_{n_2,\varepsilon})$-torsor.

Moreover, each sub-bundle $P_i$ in a decomposition of $(P,\varphi)$ 
uniquely determines the other one and comes with a self-dual isomorphism 
$\varphi_i : P_i \to P_i^*$ of the same sign as $\varphi$.
\end{proposition}

The sub-bundles $P_i$ occuring in a decomposition of $(P,\varphi)$
are characterized by the property that $\varphi(P_i)$ and
$P_i^{\perp}$ are disjoint; we then say that $P_i$ is \emph{non-degenerate}. 
A self-dual bundle will be called \emph{indecomposable}
if it admits no proper decomposition; equivalently, any proper
sub-bundle is degenerate.

\begin{remarks}\label{rem:lind}
(i) We also have the notion of $L$-\emph{indecomposability} 
from \cite{BBN05}, namely, a self-dual bundle is $L$-indecomposable 
if the associated $\PO_{n,\varepsilon}$-torsor admits no reduction of 
structure group to a proper Levi subgroup. The maximal Levi subgroups 
of $\PO_{n,\varepsilon}$ are exactly the subgroups 
$\P(\O_{n_1,\varepsilon} \times \GL_{n_2})$,
where $n_1\geq 0$, $n_2 \geq 1$, $n_1 + 2 n_2 = n$,
and $\GL_{n_2} \subset \O_{2n_2,\varepsilon}$ is the subgroup that stabilizes
a decomposition $k^{2n_2} = V_1 \oplus V_2$ with $V_1,V_2$ totally
isotropic subspaces of dimension $n_2$. Thus, a self-dual bundle is 
$L$-indecomposable if and only if it admits no proper \emph{hyperbolic} 
non-degenerate sub-bundle, where $(P,\varphi)$ is called hyperbolic 
if the bundle $P$ has a decomposition $(P_1,P_2)$ such that 
$\varphi(P_i) = P_i^{\perp}$ for $i = 1,2$. 

\smallskip

\noindent
(ii) If $P = \bP(E)$ for some vector bundle $E$ over $X$, then the
symmetric (resp. antisymmetric) morphisms $\varphi : P \to P^*$
correspond bijectively to the symmetric (resp. antisymmetric) 
non-degenerate bilinear forms $B : E \times E \to L$, where 
$L$ is a line bundle and $B$ is viewed up to multiplication by 
a regular invertible function on $X$.

Also, note that $\bP(E)$ is hyperbolic if and only if $E$ admits 
a splitting
$$
E \cong V \oplus (V^* \otimes L)
$$
for some vector bundle $V$ and some line bundle $L$; then the
bilinear form $B$ on $E$ takes values in $L$ and is given by
$$
b(v \oplus (\xi \otimes s), w \oplus (\eta \otimes t)) 
= \langle v, \eta \rangle t + \varepsilon \langle w, \xi \rangle s,
$$
where $\langle -,- \rangle$ denotes the canonical pairing on 
$V \times V^*$.  
\end{remarks}

\subsection{Structure of homogeneous self-dual bundles}
\label{subsec:str}

In this subsection, we still assume that $p \neq 2$; we denote by 
$X$ a fixed abelian variety and by $f: P \to X$ a $\bP^{n-1}$-bundle. 
We say that a self-dual bundle $(P,\varphi)$ is \emph{homogeneous}
if so is the corresponding $\PO_{n,\varepsilon}$-torsor $Z$ of 
Proposition \ref{prop:sd}. Then the bundle $P$ is easily seen 
to be homogeneous.

In view of \cite[Thm.~3.1]{Br12}, the structure of homogeneous 
self-dual bundles is described by a completely 
analogous statement to Theorem \ref{thm:str}, where $\PGL_n$ 
is replaced with $\PO_{n,\varepsilon}$. This reduces 
the classification of these bundles to that of the commutative
subgroup schemes of $\PO_{n,\varepsilon}$ up to conjugacy. Let $H$
be such a subgroup scheme, $\tH$ its preimage in $\GO_{n,\varepsilon}$
and $e: H \times H \to \bG_m$ the associated commutator pairing.
Choose a non-degenerate bilinear form $B$ on $k^n =: V$ which is 
an eigenvector of $\GO_{n,\varepsilon}$; such a form is unique up 
to scalar. We say that the pair $(\tH,B)$ is a 
\emph{self-dual theta group}, and $(V,B)$ a 
\emph{self-dual representation}.
Note that $\tH$ is equipped with a character
\begin{equation}\label{eqn:char}
\beta : \tH \longrightarrow \bG_m
\end{equation}
such that
$$
(\tx \cdot B)(v_1,v_2) = B(\tx^{-1} v_1,\tx^{-1} v_2) 
= \beta(\tx) B(v_1,v_2)
$$ 
for all $\tx \in \tH$ and $v_1,v_2 \in V$. In particular,
$\beta(t) = t^{-2}$ for all $t \in \bG_m$;
we say that $\beta$ has \emph{$\bG_m$-weight} $-2$. The existence
of such a character imposes a strong restriction on the quotient 
$H/H^{\perp} = F/F^{\perp}$ (where $F$ denotes the group of components
of $H_s$ and the orthogonals are relative to the pairing $e$):

\begin{lemma}\label{lem:se}
With the above notation, $H/H^{\perp}$ is a $2$-elementary finite group;
in particular, the homogeneous index of $P$ is a power of $2$.
Moreover, $e$ factors through a non-degenerate alternating morphism
\begin{equation}\label{eqn:se}
se: H/H^{\perp} \times H/H^{\perp} \longrightarrow \mu_2.
\end{equation}
\end{lemma}

\begin{proof}
Since $\beta(e(x,y)) = \chi(\tx \ty \tx^{-1} \ty^{-1}) = 1$ for all 
$x,y \in H$ with lifts $\tx, \ty \in \tH$,
we see that $e(2x, y) = e(x,y)^2 = 1$. Thus, $H^{\perp}$ contains 
$2H$ (the image of the multiplication by $2$ in the commutative 
group scheme $H$), i.e., $F$ is killed by $2$. Since $p \neq 2$, 
this implies the first assertion. For the second one, 
note that $e$ factors through a morphism $H \times H \to \mu_2$
and hence through a bilinear alternating morphism
(\ref{eqn:se}), which must be non-degenerate by the definition
of $H^{\perp}$.
\end{proof}

In view of this result, the statements of Proposition 
\ref{prop:theta}, Lemma \ref{lem:cent} and Proposition \ref{prop:rep}
also hold in this setting (without the assumption that $(n,p) =1$),
by the same arguments.

We now assume that $e$ is non-degenerate; equivalently, $H^{\perp}$ 
is trivial. Then we may view $H$ as a finite-dimensional vector space 
over the field $\bF_2$ with $2$ elements, and $se$ as a symplectic form 
(with values in $\bF_2$), by identifying $\bF_2$ to $\mu_2$ via 
$x \mapsto (-1)^x$. We denote by $\Sp(H) = \Aut(H,se)$ the corresponding 
symplectic group.
 
Choose a maximal totally isotropic subspace $K \subset H$. Then 
$H \cong K \oplus K^*$ and this identifies $se$ with the standard
symplectic form $\omega$ defined by
$$
\omega((x,\xi),(x',\xi')) = 
\langle x,\xi' \rangle + \langle x',\xi \rangle,
$$ 
where $\langle -, - \rangle : K \times K^* \to \bF_2$ 
denotes the canonical pairing. In particular, $\#(H) = \#(K)^2 = 2^{2r}$, 
where $r := \dim_{\bF_2}(K)$, and $\Sp(H) = \Sp_{2r}(\bF_2)$;
we say that $r$ is the \emph{rank} of $(H,e)$. Moreover,
the dual $K^*$ is identified to the character group of $K$, via the map
$\xi \longmapsto (x \mapsto (-1)^{\langle x, \xi \rangle}$).
Recall that $\tH$ is isomorphic to the Heisenberg group $\cH(K)$,
and has a unique irreducible representation of weight $1$: the 
standard representation in $\cO(K)$, of dimension $2^r$.

We now analyze the representation of $\tH$ in the space of bilinear
forms on $W$. Since $p \neq 2$, we have a decomposition
of representations 
$W^* \otimes W^* = S^2 W^* \oplus \Lambda^2 W^*$
into the symmetric and the alternating components. 
For any $x \in K$, denote by $\epsilon_x \in W^*$ the 
evaluation at $x$, i.e., $\epsilon_x(f) = f(x)$ for any $f \in W$. 
Then the $\epsilon_x$ $(x \in K)$, form a basis of $W^*$ and satisfy
$$
(t,x,\xi) \cdot \epsilon_y = 
t^{-1} \, (-1)^{\langle x + y, \xi \rangle} \, \epsilon_{x + y}.
$$
Define bilinear forms on $W$ by
$$
B_{x,\xi} := \sum_{y \in K} (-1)^{\langle y,\xi \rangle} 
\, \epsilon_y \otimes \epsilon_{x + y} \quad (x \in K, \; \xi \in K^*).
$$

\begin{lemma}\label{lem:square}
With the above notation, each $B_{x,\chi}$ is an eigenvector of 
$\tH$ with weight 
$$
\chi_{x,\xi} : (t,y,\eta) \longmapsto 
t^{-2} \, (-1)^{\langle x, \eta \rangle + \langle y,\xi \rangle}.
$$
Also, $B_{x,\xi}$ is symmetric (resp. alternating) 
if and only if $\langle x,\xi \rangle = 0$ (resp. $= 1$).
 
Moreover, the $B_{x,\xi}$ form a basis of $W^* \otimes W^*$.
\end{lemma}

\begin{proof}
The first assertion is a direct verification. It implies the
second assertion, since the $B_{x,\chi}$ have pairwise distinct
weights and their number is $\#(K)^2 = \dim(W^* \otimes W^*)$.
\end{proof}

The normalizer $N_{\GL(W)}(\tH)$ acts on $W^* \otimes W^*$; it stabilizes
$S^2 W^*$ and $\Lambda^2 W^*$, and permutes the eigenspaces of $\tH$. 
Thus, $N_{\GL(W)}(\tH)$ acts on the set of their weights,
$$
\cX := \{ \chi_{x,\xi} ~\vert~ x \in K, \; \xi \in K^* \}.
$$
Note that $\cX$ is exactly the set of characters of $\tH$ with
$\bG_m$-weight $-2$. This is an affine space with underlying vector 
space the character group of $H$, that we identify with $H$ 
via the pairing $se$. Also, $N_{\GL(W)}(\tH)$ acts on $\cX$
by affine automorphisms, and the subgroup $\tH$ of $N_{\GL(W)}(\tH)$ 
acts trivially, since $\tH$ acts on itself by conjugation. In view of 
Lemma \ref{lem:cent}, it follows that $N_{\GL(W)}(\tH)$ acts on 
$\cX$ via its quotient $\Sp(H)$; the linear part of this affine 
action is the standard action of $\Sp(H)$ on $H$. 

\begin{proposition}\label{prop:orb}
The above action of $\Sp(H)$ on $\cX$ has two orbits: 
the symmetric characters $\chi_{x,\xi}$, where 
$\langle x, \xi \rangle =0$, and the alternating characters.
In particular, $S^2 W^*$ and $\Lambda^2 W^*$ are irreducible 
representations of $N_{\GL(W)}(\tH)$.
\end{proposition}

\begin{proof}
Consider the general linear group $\GL(K) \cong \GL_r(\bF_2)$
acting naturally on $\cO(K) = W$. Then one readily checks that
this action is faithful and normalizes $\tH$; also, the resulting
homomorphism $\GL(K) \to N_{\GL(W)}(\tH)$ lifts the (injective) 
homomorphism $\GL(K) \to \Sp(H)$ associated with the natural
representation of $\GL(K)$ in $K \oplus K^*$. 
Moreover, the induced action of $\GL(K)$ on $\cX$ is given by
$\gamma \cdot \chi_{x,\xi} = \chi_{\gamma(x),\gamma(\xi)}$. 
Since the pairs $(x,\xi)$ such that $\langle x,\xi \rangle = 1$
form a unique orbit of $\GL(K)$, we see that $\Sp(H)$ acts
transitively on the alternating characters. 

On the other hand, the pairs $(x,\xi)$ such that 
$\langle x,\xi \rangle = 0$ decompose into orbits of $\GL(K)$ 
according to the (non)-vanishing of $x$ and $\xi$; this yields 
$4$ orbits if $m \geq 2$, and $3$ orbits if $m = 1$
(then the orbit with $x \neq 0 \neq \xi$ is missing). 
Note that the unique $\GL(K)$-fixed point $\chi_{0,0}$ 
(a symmetric weight) is not fixed by $\Sp(H)$: otherwise, 
the latter group would act on $\cX$ via its representation on 
$H$, and hence would act transitively on
$\cX \setminus \{\chi_{0,0} \} \cong H \setminus \{0\}$. 
But this is impossible, since $\Sp(H)$ preserves the 
symmetric weights. Also, note that $\GL(K)$ has index $2$
in its normalizer $N_{\Sp(H)}(\GL(K))$; moreover, any element
of $N_{\Sp(H)}(\GL(K))\setminus \GL(K)$ fixes $\chi_{0,0}$ and
exchanges the $\GL(K)$-orbits 
$\{ \chi_{x,0} ~\vert~ x \in K, x\neq 0 \}$
and $\{ \chi_{0,\xi} ~\vert~ \xi \in K^*, x\neq 0 \}$.

As a consequence, $\Sp(H)$ acts transitively on the symmetric
characters if $m = 1$. We now show that this property also holds 
when $m \geq 2$. In view of Lemma \ref{lem:cent}, it suffices 
to construct automorphisms $u, v \in \Aut^{\bG_m}(\tH)$ such 
that $u(\chi_{x,0}) = \chi_{x,\xi} = v(\chi_{0,\xi})$ for some
non-zero $x \in K$, $\xi \in K^*$. For this, let 
$q : K \to \bF_2$ be a quadratic form, and $\varphi : K \to K^*$
the associated alternating map, defined by
$\langle \varphi(x),y \rangle = q(x + y) + q(x) + q(y)$.
Let $u = u_q : \tH \to \tH$ be the map such that
$u(t,x,\xi) = (t \, (-1)^{q(x)}, x, \xi + \varphi(x))$.
Then one may check that $u \in \Aut^{\bG_m}(\tH)$
and $u(\chi_{x,0}) = \chi_{x,\varphi(x)}$. Since we may choose $q$ 
so that $\varphi(x) \neq 0$, this yields the desired automorphism 
$u$ (and $v$ by symmetry).
\end{proof}

By Lemma \ref{lem:square} and Proposition \ref{prop:orb},
\emph{there are exactly two isomorphism classes of self-dual 
non-degenerate theta groups of a prescribed rank}, the isomorphism
type being just the `sign'. We now construct representatives of each
class; we first consider the case of rank $1$. Then 
$H = \bF_2^2 \cong (\bZ/2\bZ)^2$ 
has a faithful homomorphism to $\PGL_2$, unique up to conjugation. 
Thus, $H$ lifts to two natural subgroups of $\GL_2$: the dihedral group 
$D \subset \O_2$, and the quaternionic group $Q \subset \Sp_2 = \SL_2$. 
Both groups are finite of order $8$; moreover, $\tH_1 := \bG_m D$ 
(resp. $\tH_0 := \bG_m Q$) is a non-degenerate theta group of rank 
$1$ equipped with a symmetric (resp. alternating) semi-invariant 
bilinear form. 

For an arbitrary rank $r$, the central product $\tH_1 \cdots \tH_1$
of $r$ copies of $\tH_1$ (the quotient of the product 
$\tH_1 \times \cdots \times \tH_1$ by the subtorus
$\{ (t_1,\ldots, t_r) ~\vert~ t_1 \cdots t_r = 1 \}$) is a 
self-dual non-degenerate theta group of rank $r$ and sign $+1$.
Similarly, the central product of $\tH_0$ with $r - 1$ copies of
$\tH_1$ is a self-dual non-degenerate theta group of rank $r$ and 
sign $-1$. 

\begin{remark}
The above description of the self-dual non-degenerate theta 
groups may also be deduced from the structure of 
\emph{extra-special $2$-groups}, i.e., of those finite groups $G$ 
such that the center $Z$ has order $2$, and $G/Z$ is $2$-elementary
(see \cite[Kap.~III, Satz~13.8]{Hu67} or 
\cite[Chap.~5, Thm.~5.2]{Go80}). Namely, by Lemma \ref{lem:se}, 
every self-dual non-degenerate theta group yields an extension
$1 \to \mu_2 \to G \to H \to 1$, where $G$ is extra-special. Yet
the approach followed here is more self-contained.
\end{remark}

Returning to an arbitrary self-dual theta group 
$(\tH \subset \GL(V),B)$, we now investigate the decomposition 
of $V$ into eigenspaces $V_{\lambda}$ of $Z(\tHs)$.
Recall from Proposition \ref{prop:rep} that 
$V_{\lambda} \cong U_{\lambda} \otimes W_{\lambda}$ as a 
representation of $\tH \cong H_u \times \tHs$, where
$W_{\lambda}$ is the standard representation of the Heisenberg group
$\tHs/\ker(\lambda)$. Also, since $B$ has weight $\beta$, we have
$B(V_{\lambda}, V_{\mu}) = \{ 0 \}$ unless $\lambda + \mu = - \beta$. 
This readily implies the following observations:

\begin{lemma}\label{lem:eigen}
{\rm (i)} As a self-dual representation, $V$ is the direct sum 
of the pairwise orthogonal subspaces $V_{\lambda}$, where 
$2 \lambda = - \beta$, and $V_{\lambda} \oplus V_{- \lambda - \beta}$, 
where $2 \lambda \neq - \beta$.   

\smallskip

\noindent
{\rm (ii)} If $2 \lambda = - \beta$, then $U_{\lambda}$ 
(resp. $W_{\lambda}$) is a self-dual representation of $H_u$ 
(resp. $\tHs/\ker(\lambda)$). Moreover, the restriction of
$B$ to $V_{\lambda}$ is the tensor product of the corresponding
bilinear forms on $U_{\lambda}$, resp. $W_{\lambda}$.

\smallskip

\noindent
{\rm (iii)} If $2 \lambda \neq - \beta$, then 
$V_{- \lambda - \beta} \cong V_{\lambda}^*(-\beta)$
as representations of $\tH$. Moreover, the restriction of
$B$ to $V_{\lambda} \oplus V_{- \lambda - \beta}$ is given by
the symmetrization or alternation of the canonical pairing
$V_{\lambda} \otimes V_{\lambda}^*( - \beta) \to k(-\beta)$.
\end{lemma}

As a direct consequence, we obtain the following analogue of 
the structure of indecomposable homogeneous bundles 
(Proposition \ref{prop:indec}):

\begin{proposition}\label{prop:indsd}
The following assertions are equivalent for a homogeneous self-dual 
bundle $(P,\varphi)$:

\smallskip

\noindent
{\rm (i)} $(P,\varphi)$ is indecomposable.

\smallskip

\noindent
{\rm (ii)} $V$ is indecomposable as a self-dual representation.

\smallskip

\noindent
{\rm (iii)} $\tHs$ is a Heisenberg group and one of the following 
cases occurs:

{\rm (I)} $V \cong U \otimes W$, where $U$ is an indecomposable 
self-dual representation of $H_u$ and $W$ is the standard irreducible 
representation of $\tHs$. Moreover, $H_s$ is $2$-elementary.

{\rm (II)} $V \cong (U \otimes W) \oplus (U^* \otimes W^*)(-\beta)$,
where $U$ is an indecomposable representation of $H_u$, $W$ is
the standard irreducible representation of $\tHs$ and 
$\beta$ is a character of $\tH_s$ of weight $-2$.
\end{proposition}  

\begin{remarks}\label{rem:Lind}
(i) In contrast to Proposition \ref{prop:indec}, there exist 
indecomposable self-dual bundles $(P,\varphi)$ such
that $\Aut_X^0(P,\varphi)$ is not unipotent. Specifically, if 
$(P,\varphi)$ is hyperbolic (typee (II) above), then the action
of $\bG_m$ on $V$ with weight spaces $U \otimes W$ of weight $1$,
and $(U^* \otimes W^*)(-\beta)$ of weight $-1$, yields a one-parameter
subgroup of bundle automorphisms of $(P,\varphi)$.

In fact, the condition that $\Aut_X^0(P,\varphi)$ is unipotent
characterizes the \emph{$L$-indecomposable} self-dual bundles. 
Also, one easily shows that the homogeneous self-dual bundle 
$(P,\varphi)$ is indecomposable if and only if the self-dual 
representation $V$ contains no non-trivial direct summand of type 
(II).

\smallskip

\noindent
(ii) If $(P,\varphi)$ is irreducible in the sense that it arises 
from a non-degenerate theta group, then $\Aut_X(P)$ is finite 
by Proposition \ref{prop:irr}; as a consequence, $\Aut_X(P,\varphi)$ 
is finite. But the converse does not hold in general, e.g., 
for homogeneous self-dual $\bP^2$-bundles associated with the 
subgroup $H$ of $\PO_3$ generated by the images of the diagonal
matrices with coefficients $\pm 1$ (then $H \cong (\bZ/2\bZ)^2$ 
and $e =0$).
Thus, the criteria for irreducibility obtained in Subsection 
\ref{subsec:st} do not extend to self-dual bundles. 
In \cite[Sec.~7.3]{BSU12}, an alternative, group-theoretical notion 
of irreducibility is introduced for homogeneous principal bundles 
under a semisimple group in characteristic $0$, and Propositions 
\ref{prop:irr} and \ref{prop:homirr} are generalized to that setting.
\end{remarks}

\end{document}